\documentstyle[12pt]{article}

\newlength{\abstractwidth}
\renewcommand{\thefootnote}{\fnsymbol{footnote}}
\renewcommand{\thanks}[1]{\footnote{#1}} 
\newcommand{\starttext}{ \setcounter{footnote}{0}
\renewcommand{\thefootnote}{\arabic{footnote}}}

\newcommand{\be}{\begin{equation}}
\newcommand{\bea}{\begin{eqnarray}}
\newcommand{\eea}{\end{eqnarray}} 
\newcommand{\ee}{\end{equation}}
 \newcommand{\<}{\langle}
\renewcommand{\>}{\rangle}
\def\ba{\begin{eqnarray}}
\def\ea{\end{eqnarray}}

\def\o{\omega}

\def\log{\,{\rm log}\,}

\def\o{\omega}

\def\o{\omega}

\def\na{\nabla}
\def\p{\partial}

\def\R{{\bf R}}

\def\ddb{{\partial\bar\partial}}

\def\na{{\nabla}}

\def\[{{\bf [}}
\def\]{{\bf ]}}

\begin{document}
\starttext \baselineskip=15pt \setcounter{footnote}{0}
\newtheorem{theorem}{Theorem}
\newtheorem{lemma}{Lemma}
\newtheorem{definition}{Definition}
\newtheorem{proposition}{Proposition}
\newtheorem{corollary}{Corollary}

\begin{center}
{\Large \bf
ANOMALY FLOWS \footnote{Work supported in part by the National Science Foundation Grants DMS-12-66033 and DMS-1605968. Key words: Green-Schwarz anomaly cancellation, quadratic terms in the curvature tensor, torsion constraints, flows, maximum principle.}}
\end{center}

\centerline{Duong H. Phong, Sebastien Picard, and Xiangwen Zhang}

\bigskip

\begin{abstract}

{\small The Anomaly flow is a flow which 
implements the Green-Schwarz anomaly cancellation mechanism
originating from superstring theory, while preserving the conformally balanced condition of Hermitian metrics. There are several versions of the flow, depending on whether the gauge field also varies, or is assumed known. A distinctive feature of Anomaly flows is that, in $m$ dimensions, the flow of the Hermitian metric has to be inferred from the flow of its $(m-1)$-th power $\o^{m-1}$. We show how this can be done explicitly, and we work out the corresponding flows for 
the torsion and the curvature tensors. The results are applied to produce criteria for the long-time existence of the flow, in the simplest case of zero slope parameter.}

\end{abstract}

\section{Introduction}
\setcounter{equation}{0}

Starting with the uniformization theorem, canonical metrics such as Hermitian-Yang-Mills and K\"ahler-Einstein metrics 
have played a major role in complex geometry. However, theoretical physics
suggests more notions of metrics which should qualify in some sense as canonical. Indeed, the classical canonical metrics are typically defined by a linear constraint in the curvature tensor. But in string theory, the key Green-Schwarz anomaly cancellation mechanism \cite{GS} for the consistency of superstring theory is an equation which involves the {\it square} of the curvature tensor. Furthermore, while the supersymmetry of the heterotic string compactified to $4$-dimensional Minkowski space-time required that the intermediate space carry a complex structure \cite{CHSW}, it allowed the corresponding Chern unitary connection to have non-vanishing torsion \cite{S}. The resulting condition is known as a Strominger system, and Calabi-Yau manifolds with their K\"ahler Ricci-flat metrics are only a special solution. What seems to emerge then is an as-yet unexplored area of non-K\"ahler geometry, where the K\"ahler condition is replaced by some specific constraint on the torsion, and the canonical metric condition is replaced by an equation on the torsion and possibly higher powers of the curvature. These equations are also novel from the point of view of the theory of partial differential equations, and it is an important problem to develop methods for their solutions.

\medskip
The goal of the present paper is to develop methods for the study of the following flow of Hermitian metrics on a $3$-dimensional complex manifold $X$, \bea
\label{anomaly1}
\p_t(\|\Omega\|_\o\o^2)&=&
i\p\bar\p\o-\alpha'({\rm Tr}Rm\wedge Rm-\Phi(t))\nonumber\\
\o(0)&=&\o_0.
\eea
Here $X$ is equipped with a nowhere vanishing $(3,0)$ holomorphic form $\Omega$, $\|\Omega\|_\o$ is the norm of $\Omega$ with respect to the Hermitian metric $\o$, defined by
\bea
\label{Omega}
\|\Omega\|_\o^2=i\Omega\wedge\bar\Omega \,\o^{-3},
\eea
and the expression $\Phi(t)$ is a given closed $(2,2)$-form
in the characteristic class $c_2(X)$, evolving with time. The expression $Rm$ is the curvature of the Chern unitary connection of $\o$, viewed as a $(1,1)$-form valued in the bundle of endomorphisms $End(T^{1,0}(X))$ of $T^{1,0}(X)$. The
initial Hermitian form $\o_0$ is required to satisfy the following {\it conformally balanced} condition
\bea
\label{balanced}
d(\|\Omega\|_{\o_0}\o_0^2)=0.
\eea

\medskip

The motivation for the flow (\ref{anomaly1}) is as follows. In \cite{S}, building on the earlier work of Candelas, Horowitz, Strominger, and Witten \cite{CHSW}, Strominger identified the following system of equations for a Hermitian metric $\o$ on $X$ and a Hermitian metric $H_{\bar\alpha\beta}$ on a holomorphic vector bundle $E\to X$, 
\bea
\label{HYM}
&&
F^{2,0}=F^{0,2}=0, \quad F\wedge\o^{2}=0\\
\label{GreenSchwarz}
&&
i\p\bar\p\o-\alpha'{\rm Tr}(Rm\wedge Rm-F\wedge F)=0\\
\label{balanced0}
&&
d^\dagger\o=i(\bar\p-\p)\log\|\Omega\|_\o,
\eea
as conditions for the product of $X$ with $4$-dimensional space-time to be a supersymmetric vacuum configuration for the heterotic string. The conditions on $F$ in the first equation above just mean that $F$ is the curvature of the Chern unitary connection of $H_{\bar\alpha\beta}$, and that $H_{\bar\alpha\beta}$ is Hermitian-Yang-Mills with respect to any metric conformal to $\o$. It is a subsequent, but basic observation of Li and Yau \cite{LY} that the third condition on $\o$ above, which is at first sight a {\it torsion constraints} condition, is equivalent to the condition that $\o$ be conformally balanced
\bea
\label{balanced1}
d(\|\Omega\|_\o\o^2)=0.
\eea
In the special case where $(X,\o)$ is a compact K\"ahler $3$-fold with $c_1(X)=0$, if we take $E=T^{1,0}(X)$, $H=\o$, then the anomaly condition is automatically satisfied. The Hermitian-Yang-Mills condition reduces to the condition that $\o$ be Ricci-flat, which can be implemented by Yau's theorem \cite{Y}. The norm $\|\Omega\|_\o$ is then constant, and the torsion constraints follow from the K\"ahler property of $\o$. Thus Calabi-Yau $3$-folds with their Ricci-flat metrics can be viewed as special solutions of the Strominger system, and they have played a major role ever since in both superstring theory and algebraic geometry \cite{CHSW}. From this point of view, it is natural to think of the pair $(\o,H)$ as a canonical metric for $(X,E)$, and if $H$ happens to be fixed for some reason, of the metric $\o$ itself as a canonical metric in non-K\"ahler geometry.

\medskip
Strominger systems are difficult to solve, and the first non-perturbative, non-K\"ahler solutions to the systems were obtained by Fu-Yau \cite{FY1, FY2}, some twenty years after Strominger's original proposal. These solutions were on toric fibrations over $K3$ surfaces constructed earlier by Goldstein and Prokushkin \cite{GP}. On such manifolds, Fu-Yau succeeded in reducing the Strominger system to a new complex Monge-Amp\`ere equation on the two-dimensional K\"ahler base, which they succeeded in solving. Higher dimensional analogues of the Fu-Yau solution were considered by the authors in \cite{PPZ2, PPZ3, PPZ4}. Geometric constructions of some special solutions of Strominger systems have been given in e.g. \cite{AGF, Fe1, Fe2, FeY, FIUV1, FIUV2, FTY, OUV}.

\medskip
A major problem at the present time is to develop analytical methods for solving the general Strominger system. Even if the curvature $F$ of the bundle metric $H$ were known and we concentrate only on the equations for $\o$, an immediate difficulty typical of non-K\"ahler geometry, is that there is no general or convenient way of parametrizing conformally balanced metrics, comparable to the parametrization of K\"ahler metrics by their potentials which was instrumental in Yau's solution of the Ricci-flat equation. It appears to be a daunting problem to have to deal with the anomaly equation and the conformally balanced equation as a system of equations. 
A way of bypassing this difficulty was suggested by the authors in \cite{PPZ1}, which is to introduce the coupled geometric flow
\bea
\label{anomaly2}
H^{-1} \, \p_t H
&=&
-\Lambda F\nonumber\\
\p_t(\|\Omega\|_\o\o^2)
&=&
i\p\bar\p\o-\alpha'{\rm Tr}(Rm\wedge Rm-F\wedge F)
\eea
with initial conditions $\o(0)=\o_0$, $H(0)=H_0$, where $H_0$ is a given metric on $E$, and $\o_0$ is a Hermitian metric on $X$ which satisfies the conformally balanced condition (\ref{balanced})\footnote{Note that there are ways for constructing individual conformally balanced metrics $\o_0$ (see e.g. \cite{TW}).}. 

\medskip

The point of the flow is that, by Chern-Weil theory, the right hand side in the second line above is always closed, and hence the condition $d(\|\Omega\|_\o\o^2)=0$ is preserved by the flow. Thus there is no need to treat the conformally balanced condition as a separate equation, and the stationary points of the flow will automatically satisfy all the equations in the Strominger system.
For fixed $\o$, the flow of the metric $H_{\bar\alpha\beta}$ is just the Donaldson heat flow \cite{D}. If the flow for $H_{\bar\alpha\beta}(t)$ is known, and if we set $\Phi(t)={\rm Tr}(F\wedge F)$, then the flow for $\o$ reduces to the flow (\ref{anomaly1}). An understanding of (\ref{anomaly1}) appears a necessary preliminary step in an understanding of (\ref{anomaly2}). The flow (\ref{anomaly2}) was called the Anomaly flow in \cite{PPZ1}, in reference to the key role played by the right hand side in the Green-Schwarz anomaly cancellation mechanism. We shall use the same generic name for all closely related flows such as (\ref{anomaly1}).

\medskip
Anomaly flows appear to be considerably more complicated than classical flows in geometry of which the Yang-Mills flow and the Ricci flow are well-known examples. A first hurdle is that the flow of metrics $\o(t)$ has to be deduced from the flow of $(2,2)$-forms $\|\Omega\|_{\o}\,\o^2$. Now the existence in dimension $m$ of an $(m-1)$-th root of a positive $(m-1,m-1)$-form has been shown by Michelsohn \cite{M}, and this passage back and forth between positive $(1,1)$-forms and $(m-1,m-1)$-forms
has played a major role e.g. in works of Popovici \cite{Po} and in the recent proof by Szekelyhidi, Tosatti, and Weinkove \cite{STW, TW} of the existence of Gauduchon metrics with prescribed volume form. However, it does not appear possible to use the formalism in these works to deduce the flow of the curvature tensor of $\o$ from the flow of $\o^{m-1}$. This is one of the main goals of the present paper. 
What we do is to produce a seemingly new formula for the square root of a $(2,2)$-form, or equivalently, for the Hodge $\star$ operator, without using the antisymmetric symbol $\varepsilon$. 
With such a formula, and using the very specific torsion constraints resulting from the conformally balanced condition, we obtain the following completely explicit expression for the Anomaly flow:

\begin{theorem}
\label{anomaly4}
If the initial metric $\o_0$ is conformally balanced, then
the Anomaly flow (\ref{anomaly1}) can also be expressed as
\bea
\p_t g_{\bar pq}
=
{1\over 2\|\Omega\|_\o}
\bigg[-\tilde R_{\bar pq}+g^{\alpha\bar\beta}g^{s\bar r}
T_{\bar\beta s q}\bar T_{\alpha \bar r\bar p}
-
\alpha'
g^{s\bar r}(R_{[\bar ps}{}^\alpha{}_\beta R_{\bar rq]}{}^\beta{}_\alpha
-\Phi_{\bar ps\bar rq})\bigg].
\eea
where $\tilde R_{\bar kj}$ is the Ricci tensor and $T_{\bar k i j}$ is the torsion tensor, as defined in (\ref{Ricci-tensor}) and (\ref{torsion-tensor}) below. The brackets $[\, ,\,]$ denote anti-symmetrization separately in each of the two sets of barred and unbarred indices. 
\end{theorem}

The above theorem shows that the Anomaly flow can be viewed as generalization of the Ricci flow, with higher order corrections in the curvature tensor proportional to $\alpha'$. Indeed, the terms $\tilde R_{\bar pq}-g^{\alpha\bar\beta}g^{s\bar r}
T_{\bar\beta s q}\bar T_{\alpha \bar r\bar p}$ reduce to the Ricci curvature $R_{\bar pq}$ (see the definition in (\ref{Ricci-tensor})) if the torsion vanishes, and the terms with coefficient $\alpha'$ are the higher order corrections. It is remarkable that this analogy with the Ricci flow is due not to an attempt to generalize the Ricci flow, but rather to the combination of the Green-Schwarz cancellation mechanism, more specifically the de Kalb-Ramond field $i\p\bar\p\o$, with the torsion constraints equivalent to the conformally balanced condition.
Once the formulation of the flow provided by Theorem \ref{anomaly4} is available, it is straightforward to derive the flows of the torsion and curvature tensors. The full results are given in Theorems \ref{curvatureflow} and \ref{torsionflow} below. Here, we note only that they reinforce the same analogy with the Ricci flow. For example, the diffusion operator in the flow for the Ricci curvature is given by
\bea
\p_tR_{\bar kj}
=
{1\over 2\|\Omega\|_\o}\big(\Delta R_{\bar kj}
+
2\alpha' g^{\lambda\bar\mu}g^{s\bar r}R_{[\bar r\lambda}{}^\beta{}_\alpha
\na_s\na_{\bar\mu]} R_{\bar kj}{}^\alpha{}_\beta)
+\cdots
\eea
Up to the factor $(2\|\Omega\|_\o)^{-1}$, it coincides with the diffusion operator $\Delta$ for the Ricci curvature in the Ricci flow, up to a higher order correction in the curvature which is proportional to $\alpha'$.

\medskip
The formulation of the Anomaly flow provided by Theorem \ref{anomaly4} makes it more amenable to existing techniques for flows, and indeed many flows of metrics with torsion have been studied in the literature (e.g. \cite{Gi, LW, SW, ST, TW1} and others). However, the Anomaly flow still involves a combination of novel features such as the particular torsion constraints, the presence of the factor $\|\Omega\|_\o$ (which is quite important in string theory as it originates from the dilaton field), and especially the presence of the quadratic terms in the curvature tensor. All this makes a general solution only a remote possibility at this time. Thus we focus on two important special cases. The first case is the Anomaly flow restricted to the Fu-Yau ansatz for solutions of the Strominger system on toric fibrations over Ricci-flat K\"ahler surfaces. We can show that the Anomaly flow converges in this case, and thus gives another proof of the existence theorem of Fu-Yau. But because of its length and complexity, the full argument will be presented in a companion paper \cite{PPZ5}  to the present one. The second case is when $\alpha'=0$. 
In this case, our main results are as follows. 

\begin{theorem}\label{Higherorderestimate}
Assume that $\alpha'=0$. Suppose that $A>0$ and $\o(t)$ is a solution to the Anomaly flow (\ref{anomalyK}) below, with $t\in [0, {1\over A}]$. Then, for all $k\in \mathbf{N}$, there exists a constant $C_k$ depending on a uniform lower bound of $\|\Omega\|_{\o}$ such that, if
\bea\label{assumptiononcurvature}
|Rm|_{\o} + |D T|_{\o} + |T|^2_{\o} \leq A,    \ \ \ \textit{ for all } z\in M \textit{ and } t\in [0, {1\over A}],
\eea
then,
\bea\label{CTestimate}
|D^k Rm(z, t)|_{\o} \leq {C_k A\over t^{k/2}}, \ \ \ |D^{k+1} T(z, t)|_{\o} \leq  {C_k A\over t^{k/2}}
\eea
for all $z\in M$ and $t\in (0, {1\over A}]$.
\end{theorem}

The estimates given in the above theorem can be viewed as Shi-type derivative estimates for the curvature tensor and torsion tensor along the Anomaly flow (\ref{anomalyK}). With this theorem, we can provide a criterion for the long-time existence of the Anomaly flow:

\begin{theorem}
\label{longtimeexistence}
Assume that $\alpha'=0$, and that the Anomaly flow (\ref{anomalyK}) exists on an interval $[0,T)$ for some $T>0$. If ${\rm inf}_{t\in [0,T)}\|\Omega\|_\o>0$ (or equivalently $\o^3(t)\leq C\,\o^3(0)$), and if
\bea
{\rm sup}_{X\times [0,T)}
(|Rm|_{\o}^2+|DT|_{\o}^2+|T|_{\o}^4)<\infty
\eea
then the flow can be continued to an interval $[0,T+\epsilon)$ for some $\epsilon>0$. In particular, the flow exists for all time, unless there is a time $T>0$ and a sequence $(z_j,t_j)$, with $t_j\to T$, with either $\|\Omega(z_j,t_j)\|_\o\to 0$,
or 
\bea
(|Rm|_{\o}^2+|DT|_{\o}^2+|T|_{\o}^4)(z_j,t_j)\to \infty.
\eea
\end{theorem} 

\medskip
The paper is organized as follows. In \S 2, we begin by providing an effective way for recapturing the form $\p_t\o$ from the form $\p_t(\|\Omega\|_\o\o^2)$. We then discuss the torsion constraints in the Strominger system, and in particular, how they result in two different notions of Ricci curvature, but a single notion of scalar curvature. We can then prove Theorem \ref{anomaly4}. With Theorem \ref{anomaly4}, it is straightforward to derive the flows of the curvature and of the torsion. In \S3, we give the proof of Theorem \ref{Higherorderestimate}. This proof is analogous to the proof for the classical flows, but it is more complicated here due to the non-vanishing torsion and the expression $\|\Omega\|_\o$. Once we have Theorem \ref{Higherorderestimate}, it is easy to prove Theorem \ref{longtimeexistence}. Finally, we provide a list of conventions in the appendices, together with some basic identities of Hermitian geometry.

\


\section{The flows of the metric, torsion and curvature}
\setcounter{equation}{0}

The first task in the study of a geometric flow is to derive the flows of the curvature tensor, and in the case of non Levi-Civita connections, of the torsion tensor. 
In the case of Anomaly flows, this task is complicated by the fact that the flow is defined as a flow of the $(2,2)$-form $\|\Omega\|_\o \,\o^2$, and that the flow of $\o$ itself has to be recaptured from there.

\medskip

\subsection{The equation $\varphi\wedge \o^{m-2}=\Phi$ and the Hodge $\star$ operator}

Since $\p_t\o^2=2\p_t\o\wedge\o$, the flow of $\o$ can be recovered from the flow of $\o^2$ if we can solve explicitly equations of the form $\varphi\wedge\o= \Phi$ for a given $\Phi$. We begin by doing this, in general dimension $m$ instead of just $m=3$, as the resulting formulas for the solution as well as the Hodge $\star$ operator may be of independent interest.

Let $\o=ig_{\bar kj}dz^j\wedge d\bar z^k$ and $\eta$ be a $(p,q)$-form. We define its components
$\eta_{\bar k_1\cdots\bar k_q j_1\cdots j_p}$ by
\bea
\eta=
{1\over p!q!}
\sum \eta_{\bar k_1\cdots\bar k_q j_1\cdots j_p}\,
dz^{j_p}\wedge\cdots\wedge dz^{j_1}\wedge
d\bar z^{k_q}\wedge\cdots\wedge d\bar z^{k_1}.
\eea

\begin{lemma}
\label{factorization}
Let $\Phi$ be a $(m-1,m-1)$ form on a Hermitian manifold $(X,\o)$ of dimension $m$. Then the equation
\bea
\varphi\wedge  \o^{m-2}=\Phi
\eea
admits a unique solution, given by
\bea
\label{star}
\varphi_{\bar jk}
=
{1\over \alpha_m}
\bigg\{
i^{-(m-2)}\, \prod_{p=1}^{m-2}g^{k_p\bar j_p}\Phi_{\bar jk\bar j_1k_1\cdots\bar j_{m-2}k_{m-2}}
-
{\beta_m \over (m-1)!^2} ({\rm Tr}\,\Phi)\,i\,  g_{\bar jk}\bigg\}
\eea
where $\alpha_m$ and $\beta_m$ are universal constants, depending only on the dimension $m$, given by
\bea
\alpha_m=(m-1)!(m-2)!(m-1-{m^2\over 6}),
\qquad
\beta_m={m!(m-2)!\over 6}.
\eea
and ${\rm Tr}\,\Theta$ for a $(p,p)$-form $\Theta$ is defined by
\bea\label{defineTr}
{\rm Tr}\,\Theta=\<\Theta, \o^p\>=
i^{-p}\, \prod_{\ell=1}^p g^{k_\ell\bar j_\ell}\,\Theta_{\bar j_1k_1\cdots\bar j_pk_p}.
\eea
The traces of $\varphi$ and $\Phi$ are related by
\bea
{\rm Tr}\varphi={1\over (m-1)!^2}{\rm Tr}\,\Phi.
\eea
\end{lemma}

\medskip
\noindent
{\it Proof.} In components, the equation $\varphi\wedge \o^{m-2}=\Phi$ can be expressed as
\bea
i^{m-2}\,\varphi_{\{\bar jk}g_{\bar j_1k_1}\cdots g_{\bar j_{m-2}k_{m-2}\}}
=
\Phi_{\bar jk\bar j_1k_1\cdots\bar j_{m-2}k_{m-2}}
\eea
where the bracket $\{,\}$ denote antisymmetrization of all the barred indices as well as of all the unbarred indices. We contract both sides, getting
\bea
i^{m-2}\, \prod_{p=1}^{m-2}g^{k_p\bar j_p}\varphi_{\{\bar jk}g_{\bar j_1k_1}\cdots g_{\bar j_{m-2}k_{m-2}\}}
=
\prod_{p=1}^{m-2}g^{k_p\bar j_p}\Phi_{\bar jk\bar j_1k_1\cdots\bar j_{m-2}k_{m-2}}.
\eea
We expand the left-hand side by writing down all the terms arising from antisymmetrization of the sub-indices. Carrying out the contractions with
$\prod_{p=1}^{m-2}g^{k_p\bar j_p}$, it is easy to verify that each term is a constant multiple of $\varphi_{\bar jk}$ or of (${\rm Tr}\,\varphi)\,g_{\bar kj}$. This shows that we have a relation of the form
\bea
\label{star1}
i^{m-2}\, \alpha_m \varphi_{\bar kj}+i^{m-1}\, \beta_m ({\rm Tr}\,\varphi)\,g_{\bar kj}
=
\prod_{p=1}^{m-2}g^{k_p\bar j_p}\Phi_{\bar jk\bar j_1k_1\cdots\bar j_{m-2}k_{m-2}}.
\eea
Next, we have $\varphi\wedge \o^{m-1}=\Phi\wedge\o$,
which implies
\bea
\<\varphi,\star\o^{m-1}\>=\<\Phi,\star\o\>
\eea
Recalling that $\star\,\o^{m-1}=(m-1)!\, \o$ and $\star\,\o={1\over (m-1)!}\, \o^{m-1}$, we obtain
\bea
{\rm Tr}\,\varphi=\<\varphi, \o\>={1\over (m-1)!^2}\,{\rm Tr}\,\Phi.
\eea
This establishes the form (\ref{star}).

\smallskip
It is easy to see that $\alpha_m\not=0$, otherwise we obtain a relation between $g_{\bar kj}$ and $\Phi$ that cannot hold for an arbitrary $(m-1,m-1)$-form $\Phi$. To determine the precise values of $\alpha_m$ and $\beta_m$, we proceed as follows.

\smallskip
First, contracting (\ref{star1}) with respect to $g^{k\bar j}$ and using the definition of ${\rm Tr}\Theta$ in (\ref{defineTr}), we have the following relation between $\alpha_m$ and $\beta_m$,
\bea
i^{m-1}\,\alpha_m\, {\rm Tr}\,\varphi +i^{m-1}\, \beta_m m\, {\rm Tr}\,\varphi
=
i^{m-1}\,{\rm Tr}\,\Phi=i^{m-1}\,(m-1)!^2\, {\rm Tr}\,\varphi,
\eea
and hence
\bea
\label{star0}
\alpha_m+\beta_m m=(m-1)!^2.
\eea
Thus it remains only to determine $\beta_m$. We note that the only permutation of indices which can produce a multiple of $({\rm Tr}\,\varphi)\,g_{\bar kj}$ is of the form, e.g., 
\bea
\label{star2}
\varphi_{\bar j_1k_1}g_{\bar jk}\,g_{\bar j_2k_2}\cdots g_{\bar j_{m-2}k_{m-2}}
\eea
which produces, upon contraction with $\prod_{p=1}^{m-2}g^{k_p\bar j_p}$
and antisymmetrization in $j_2,\cdots,j_{n-2}$ and $k_2,\cdots,k_{n-2}$,
\bea
({\rm Tr}\,\varphi) g_{\bar kj}\,\<\o^{m-3},\o^{m-3}\>
=
({\rm Tr}\,\varphi) g_{\bar kj}\|\o^{m-3}\|^2.
\eea
We can compute $\|\o^{m-3}\|^2$ as follows
\bea
\o^{m-3}
=
(m-3)!\sum_{j<k<\ell} (ie^j\wedge \overline{e^j})^{v}\cdots
(ie^k\wedge \overline{e^k})^{v}\wedge
(ie^\ell\wedge \overline{e^\ell})^{v}
\eea
Since the sum in the right hand side consists of exactly ${1\over 3!}m(m-1)(m-2)$ terms, we find
\bea
\|\o^{m-3}\|^2 =(m-3)!^2{m(m-1)(m-2)\over 6}
=
{m!(m-3)!\over 6}.
\eea
But there are $m-2$ terms of the form (\ref{star2}), corresponding to the indices $\bar j_1k_1$ taking successively all values to $\bar j_{m-2}k_{m-2}$. Thus we obtain
\bea
\beta_m={m!(m-3)!\over 6}(m-2)={m!(m-2)!\over 6}
\eea
establishing our claim for $\beta_m$. The claim for $\alpha_m$ then follows from the relation (\ref{star0}).

\

Although the previous lemma suffices for our purpose, it is useful for future considerations to point out that it gives in effect an explicit expression for the Hodge $\star$ operator without the $\epsilon$ symbol. This can be seen by comparing it with the following lemma which solves the same equation, but using the Hodge $\star$ operator (and which can also be derived from Proposition 1.2.31 in \cite{Huy}):

\begin{lemma}
Let $(X,\o)$ be a Hermitian manifold of complex dimension $m\geq 2$. Consider the following equation, for a given $(m-1,m-1)$-form $\Phi$,
\bea
\label{factorization_m}
\psi\wedge \o^{m-2}=\Phi.
\eea
Then the equation admits a unique solution, given by
\bea
\psi=-{1\over (m-2)!}\star\Phi+{\<\Phi,\o^{m-1}\>\over (m-1)!^2}\, \o
\eea
\end{lemma}

\smallskip
\noindent
{\it Proof.} First observe that if $\psi_0$ is a $(1,1)$-form with $\<\psi_0,\o\>=0$, then
\bea
\star(\psi_0\wedge\o^{m-2})=-(m-2)!\psi_0
\eea
as can be verified by working out $\psi_0\wedge\o^{m-2}$. This is equivalent to saying that, if $\Phi_0$ is a $(m-1,m-1)$ form with $\<\Phi_0,\o^{m-1}\>=0$, then $\psi_0=-{1\over (m-1)!}\star\Phi_0$ is the unique solution of the equation
$\psi_0\wedge\o^{m-2}=\Phi_0$. Next, for general $\Phi$, we write
\bea
\Phi={\<\Phi,\o^{m-1}\>\over\|\o^{m-1}\|^2}\o^{m-1}+\Phi_0
\eea
so that $\<\Phi_0,\o^{m-1}\>=0$. In view of the previous observation, the $(1,1)$-form
\bea
\psi&=&-{1\over (m-2)!}\star\Phi_0+{\<\Phi,\o^{m-1}\>\over\|\o^{m-1}\|^2}\o
\nonumber\\
&=&
-{1\over (m-2)!}\star(\Phi-{\<\Phi,\o^{m-1}\>\over\|\o^{m-1}\|^2}\o^{m-1})
+
{\<\Phi,\o^{m-1}\>\over\|\o^{m-1}\|^2}\o
\nonumber\\
&=&
-{1\over (m-2)!}\star\Phi+
m{\<\Phi,\o^{m-1}\>\over\|\o^{m-1}\|^2}\o
\eea
is a solution of the equation (\ref{factorization_m}). Here we used the fact that $\star\o^{m-1}=(m-1)!\o$. Since $\|\o^{m-1}\|^2=m!(m-1)!$, we obtain the desired formula. Q.E.D.

\

Comparing the previous two lemmas gives the following formula for the Hodge $\star$ operator on $(m-1,m-1)$ forms, on an arbitrary Hermitian $m$-fold $(X,\o)$,

\begin{lemma}
Let $(X,\o)$ be a Hermitian manifold of complex dimension $m\geq 2$. Then for any $(m-1,m-1)$-form $\Phi$, we have
\bea
(\star\Phi)_{\bar j k}
=
{1\over (m-1)!(m^2-6m+6)}
\bigg\{
6\, i^{-(m-2)}\, \prod_{p=1}^{m-2}g^{k_p\bar j_p}\Phi_{\bar jk\bar j_1k_1\cdots\bar j_{m-2}k_{m-2}}
+(m-6)({\rm Tr}\Phi)\, i \, g_{\bar j k}\bigg\}
\nonumber
\eea
\end{lemma}

\medskip
\noindent
{\it Proof.} We equate $\psi= \varphi$, in the notation of the previous two lemmas. Thus
\bea
&&-{1\over (m-2)!}(\star\Phi)_{\bar jk}
+
{\<\Phi,\o^{m-1}\>\over (m-1)!^2}ig_{\bar jk}\nonumber\\
&=&
\bigg\{
{i^{-m+2} \over\alpha_m}\prod_{p=1}^{m-2}g^{k_p\bar j_p}\Phi_{\bar jk\bar j_1k_1\cdots\bar j_{m-2}k_{m-2}}
-
{\beta_m\over\alpha_m}{i\over (m-1)!^2}({\rm Tr}\,\Phi)g_{\bar jk}\bigg\}.
\eea
Since ${\rm Tr}\,\Phi =\<\Phi,\o^{m-1}\>$, we obtain
\bea
-{1\over (m-2)!}\star\Phi
=
{i^{-m+2}\over\alpha_m}\prod_{p=1}^{m-2}g^{k_p\bar j_p}\Phi_{\bar jk\bar j_1k_1\cdots\bar j_{m-2}k_{m-2}}
-
(1+{\beta_m\over\alpha_m}){{\rm Tr}\Phi \over (m-1)!^2}\, i \, g_{\bar j k}
\eea
Working out the coefficient $1+\beta_m/\alpha_m$, we obtain the formula
\bea
{1\over (m-2)!}
\star\Phi
=-{i^{-m+2}\over\alpha_m}
\prod_{p=1}^{m-2}g^{k_p\bar j_p}\Phi_{\bar jk\bar j_1k_1\cdots\bar j_{m-2}k_{m-2}}
+
{(m-1)(m-6)\over m^2-6m+6}{{\rm Tr}\Phi \over (m-1)!^2}\, i \, g_{\bar j k}.
\nonumber
\eea
This can be rewritten in turn in the form given in the lemma. Q.E.D.

\

\subsection{Torsion and curvature for conformally balanced metrics}

Next, we examine more carefully the implications for the torsion and curvature condition of conformally balanced metrics. Let $\o=ig_{\bar kj}dz^j\wedge d\bar z^k$ be a Hermitian metric, viewed as a positive $(1,1)$-form. We define its torsion tensor $T$ and $\bar T$ by
\bea
T=i\p\o,\qquad\bar T=-i\bar\p\o
\eea
which are respectively $(2,1)$ and $(1,2)$ forms. Following the conventions for $(p,q)$-forms given in the appendix, we define the coefficients $T_{\bar kjm}$ and $\bar T_{j\bar p\bar q}$ by
\bea
T={1\over 2}T_{\bar kjm}dz^m\wedge dz^j\wedge d\bar z^k,
\quad
\bar T={1\over 2}\bar T_{k\bar j\bar m}d\bar z^m\wedge d\bar z^j\wedge dz^k,
\eea
and thus
\bea\label{torsion-tensor}
T_{\bar kjm}=\p_jg_{\bar km}-\p_mg_{\bar kj},
\quad
\bar T_{k\bar j\bar m}=\p_{\bar j}g_{\bar mk}-\p_{\bar m}g_{\bar jk}.
\eea
For our later use, it is convenient to introduce
\bea
T_m = g^{j \bar k} T_{\bar k j m},
\quad
\bar T_{\bar m} = g^{\bar j k} \bar T_{k \bar j \bar m}.
\eea

\medskip

As noted earlier, the equivalence between (\ref{balanced0}) and (\ref{balanced1}) had been pointed out by Li and Yau \cite{LY}. For the convenience of the reader, we provide here the proof, stated in a somewhat more general form.

\begin{lemma}
\label{cb}
Let $(X,\o)$ be a $m$ dimensional Hermitian manifold equipped with a nowhere vanishing holomorphic $(m,0)$-form $\Omega$. Then the following conditions are equivalent:

{\rm (i)} The metric $\o$ satisfies the conformally balanced condition 
$d(\|\Omega\|_\o^a \o^{m-1})=0$ for some constant $a$;

{\rm (ii)} $d^\dagger\o=i(\bar\p-\p)\log\|\Omega\|^a$;

{\rm (iii)} $T_q=\p_q\log\|\Omega\|_\o^a$, $\bar T_{\bar q}=\p_{\bar q}\log\|\Omega\|_\o^a$.
\end{lemma}

\smallskip
\noindent
{\it Proof.} The conformally balanced condition can be written as
\bea
\p\log \|\Omega\|_\o^a\wedge \o^{m-1}+(m-1)\p\o\wedge\o^{m-2}=0.
\eea
Now just as $i\theta \wedge \o^{m-1}=-i(g^{j\bar k}\theta_{\bar kj}){\o^m\over m}$ for any $(1,1)$-form $\theta$, it is easy to verify that
\bea
T\wedge \o^{m-2}=
-i(g^{j\bar k}T_{\bar kjm} dz^m)\wedge {\o^{m-1}
\over (m-1)}
\eea
for any $(2,1)$-form $T$. Substituting $T=i\p\o$, and using the previous equation gives
\bea
(\p\log\|\Omega\|_\o^a-T_pdz^p)\wedge \o^{m-1}=0,
\eea
which implies $\p\log\|\Omega\|_\o^a-T_pdz^p=0$ and proves the equivalence between (i) and (iii). Finally, the equivalence between (ii) and (iii) follows at once from the expressions of the adjoints of $\p$ and $\bar\p$ on $(1,1)$-forms for a Hermitian metric
\bea
(\bar\p^\dagger \Phi)_q=
g^{k\bar p}(\na_k\Phi_{\bar pq}-T_k\,\Phi_{\bar pq}),
\qquad
(\p^\dagger\Phi)_{\bar q}
=
-g^{p\bar j}(\na_{\bar j}\Phi_{\bar qp}-\bar T_{\bar j}\Phi_{\bar qp})
\eea
In particular, when $\Phi=\o$, $\Phi_{\bar pq}=ig_{\bar pq}$, we obtain
\bea
(\bar\p^\dagger\o)_q
=
-iT_q,
\qquad
(\p^\dagger\o)_{\bar q}=i\bar T_{\bar q}.
\eea
This implies the equivalence between (ii) and (iii). Q.E.D.

\

We turn to the notion of Ricci curvature for conformally balanced metrics.
Although we have a single notion of Riemann curvature tensor, 
\bea
R_{\bar kj}{}^p{}_q=-\p_{\bar k}(g^{p\bar \ell}\p_jg_{\bar\ell q}),
\qquad
Rm = R_{\bar kj}{}^p{}_q dz^j\bar dz^k
\in \Lambda^{1,1}\otimes End(T^{1,0}(X)),
\eea
the lack of the standard symmetries for Levi-Civita connections leads to 4 different notions of Ricci curvature, defined as follows
\bea\label{Ricci-tensor}
&&
R_{\bar kj}=R_{\bar kj}{}^p{}_p,\qquad \tilde R_{\bar kj}=R^p{}_p{}_{\bar kj},
\qquad
R_{\bar kj}'=R_{\bar k}{}^p{}_p{}_j,
\qquad
R_{\bar kj}''=R^p{}_{j\bar k}{}_p.
\eea
Corresponding to these 4 notions of Ricci curvature are 4 notions of scalar curvature
\bea
\label{ricci}
R=g^{j\bar k}R_{\bar kj},
\quad
\tilde R=g^{j\bar k}\tilde R_{\bar kj},
\quad
R'=g^{j\bar k}R_{\bar kj}',
\quad
R''=g^{j\bar k}R_{\bar kj}''.
\eea

With the help from the torsion constraints, we have some nice relation between these different notions of Ricci curvature and scalar curvature. The following lemma is essential for our subsequent calculations:

\begin{lemma}
\label{torsion/curvature}
Assume that $\o$ is a conformally balanced metric on an $m$-fold $X$, in the sense that the equivalent conditions in Lemma \ref{cb} are satisfied. Then

{\rm (i)} $\na_{\bar k}T_j=\na_j\bar T_{\bar k}={a\over 2}R_{\bar kj}$.

{\rm (ii)} $R_{\bar kj}'=R_{\bar kj}''=(1-{a\over 2})R_{\bar kj}$.

{\rm (iii)} $\tilde R_{\bar kj}=(1-{a\over 2})R_{\bar kj}+\na^mT_{\bar kjm}$.

{\rm (iv)} $R=\tilde R$ and $R'=R''=(1-{a\over 2})R$.

\end{lemma}

\smallskip
\noindent
{\it Proof.} By definition, $R_{\bar kj}=-\p_j\p_{\bar k}\log \o^m=
\p_j\p_{\bar k}\log \|\Omega\|_\o^2$, so (i) follows from (iii) in Lemma \ref{cb}. Next,
\bea
R_{\bar kj}'
=R_{\bar kp}{}^p{}_j=R_{\bar kj}+\na_{\bar k}T^m{}_{jm}
=R_{\bar kj}-\na_{\bar k}T_j
=(1-{a\over 2})R_{\bar kj},
\eea
which proves a first part of (ii). Similarly,
\bea
R_{\bar kj}''=R_{\bar kj}{}^p{}_p+\na_j\bar T^{\bar q}{}_{\bar k\bar q}
=(1-{a\over 2})R_{\bar kj}
\eea
which completes the proof of (ii). Next,
\bea
\tilde R_{\bar kj}
=R_{\bar kj}+\na_j\bar T^{\bar p}{}_{\bar k\bar p}
+
\na_{\bar p}T_{\bar kjm} g^{m\bar p}
=
(1-{a\over 2})R_{\bar kj}+\na^mT_{\bar kjm}
\eea
which proves (iii). Contracting with $g_{\bar kj}$, we obtain (iv). Q.E.D. 

Note that the identities $R=\tilde R$ and $R'=R''$ actually hold for any Hermitian metric.

\

\subsection{Flow of the metric $\o$ and proof of Theorem \ref{anomaly4}}

We can now come back to the derivation of the flow for the metric $\o$ in the Anomaly flow and prove Theorem \ref{anomaly4}. In the following, we let the dimension $m$ of the manifold $X$ be then $3$, and take $a=1$ in Lemma \ref{cb}. 

\smallskip

\par
It is convenient to denote the right hand side of the Anomaly flow by $\Psi$,
\bea
\label{Psi}
\Psi=i\p\bar\p\o-\alpha' \,{\rm Tr}(Rm\wedge Rm-\Phi(t))
\eea
which is then a $(2,2)$-form. As usual, we denote its coefficients by $\Psi_{\bar ps\bar r q}$, and also introduce the notation $\Psi_{\bar pq}$, which can be viewed as the coefficients of a $(1,1)$-form,
\bea
\label{Psi0}
\Psi=
{1\over (2!)^2}\sum
\Psi_{\bar ps\bar rq}dz^q\wedge d\bar z^r\wedge dz^s\wedge d\bar z^p,
\qquad
\Psi_{\bar pq}=g^{s\bar r}\Psi_{\bar ps\bar rq}.
\eea

We rewrite the Anomaly flow (\ref{anomaly1}) as
\bea
\label{anomaly3}
(\p_t\log\|\Omega\|_\o \o+2\p_t\o)\wedge\o={1\over\|\Omega\|_\o} \Psi.
\eea

We apply the second statement in Lemma \ref{factorization}. Since 
\bea
\p_t\log\|\Omega\|_\o=-{1\over 2}\p_t\log({\rm det}\,\o)=-{1\over 2}{\rm Tr}(\p_t\o)
\eea
we find, in dimension $m=3$,
\bea
{\rm Tr}(\p_t\o)={1\over 2\|\Omega\|_\o}{\rm Tr}\Psi.
\eea
This gives us the flow of the volume form $\o^3$. Returning once again to the flow (\ref{anomaly3}) and applying the first statement in Lemma \ref{factorization}, we find
\bea \label{anomaly03}
\p_t g_{\bar pq}={1\over 2\|\Omega\|_\o} g^{s\bar r}\Psi_{\bar ps\bar rq}
=
{1\over 2\|\Omega\|_\o}\Psi_{\bar pq}.
\eea

It remains to work out the components $\Psi_{\bar ps\bar rq}$ more explicitly. The first term is
\bea
i\ddb\o
=
{1\over 2^2}\bigg\{\p_{\bar k}(\p_jg_{\bar\ell m}-\p_mg_{\bar\ell j})
-
\p_{\bar\ell}(\p_jg_{\bar km}-\p_mg_{\bar \ell j})\bigg\}
d\bar z^\ell\wedge dz^j\wedge dz^m\wedge d\bar z^k
\eea
and hence
\bea
(i\ddb\o)_{\bar kj\bar\ell m}
=
\p_{\bar\ell}(\p_jg_{\bar km}-\p_mg_{\bar kj})
-
\p_{\bar k}(\p_jg_{\bar \ell m}-\p_mg_{\bar\ell j}).
\eea

On the other hand, the Riemann curvature tensor is given by
\bea
R_{\bar kj}{}^\ell{}_m
=
-\p_{\bar k}(g^{\ell\bar p}\p_jg_{\bar p m})
=
-g^{\ell\bar p}\p_{\bar k}\p_jg_{\bar p m}
+
g^{\ell\bar r}\p_{\bar k}g_{\bar rs}g^{s\bar q}\p_jg_{\bar qm},
\eea
or, equivalently,
\bea
R_{\bar kj\bar \ell m}
=
-
\p_{\bar k}\p_jg_{\bar \ell m}
+
\p_{\bar k}g_{\bar \ell s}g^{s\bar r}\p_jg_{\bar r m}.
\eea
Thus we obtain
\bea \label{localomega}
(i\ddb\o)_{\bar kj\bar\ell m}
&=&
R_{\bar kj\bar\ell m}-R_{\bar km\bar\ell j}+R_{\bar\ell m\bar kj}-R_{\bar\ell j\bar km}
+g^{s\bar r}\,T_{\bar r mj}\bar T_{s\bar k\bar\ell}.
\eea
Applying Lemma \ref{torsion/curvature} on the torsion and Ricci curvatures of conformally balanced metrics gives
\bea \label{iddb-omega-identity}
g^{m\bar\ell}
(i\ddb\o)_{\bar kj\bar\ell m}
=
\tilde R_{\bar kj}-g^{s\bar r}g^{m\bar\ell}T_{\bar r mj}\bar T_{s\bar \ell\bar k}.
\eea

We collect the resulting formulas in a lemma:

\begin{lemma}
\label{Psi1}
Let the $(2,2)$-form $\Psi$ be defined by (\ref{Psi}) and its components $\Psi_{\bar ps\bar rq}$, $\Psi_{\bar pq}$ by (\ref{Psi0}). Then 
\bea
\label{Psi2}
&&
\Psi_{\bar km\bar \ell j}
=
R_{\bar km\bar\ell j}-R_{\bar kj\bar\ell m}+R_{\bar\ell j\bar km}-R_{\bar\ell m\bar kj}
+g^{s\bar r}\,T_{\bar r jm}\bar T_{s\bar k\bar\ell}
-
\alpha'(R_{[\bar km}{}^\alpha{}_\beta R_{\bar \ell j]}{}^\beta{}_\alpha-
\Phi_{\bar km\bar \ell j})
\nonumber\\
&&
\Psi_{\bar k j}
=
-\tilde R_{\bar kj}+(T\bar T)_{\bar kj}
-\alpha' g^{m\bar\ell}
(R_{[\bar k m}{}^\alpha{}_\beta R_{\bar \ell  j]}{}^\beta{}_\alpha-\Phi_{\bar km\bar\ell j})
\eea
where the brackets $[\,,\,]$ denote anti-symmetrization separately in each of the two sets of barred and unbarred indices and $(T\bar T)_{\bar kj}:= g^{s\bar r}g^{m\bar\ell}T_{\bar r mj}\bar T_{s\bar \ell\bar k}$.
\end{lemma}

Combining the formula (\ref{anomaly03}) for the Anomaly flow, and using the fact that the flow preserves the conformally balanced condition, we obtain Theorem \ref{anomaly4}.

\

\subsection{Flow of the curvature tensor}

The general formula for the flow of the curvature tensor of Chern unitary connections under a flow of metrics is the following
\bea
\p_t R_{\bar kj}{}^\mu{}_\nu
=
-\na_{\bar k}\na_j (g^{\mu\bar\gamma}\dot g_{\bar\gamma\nu})
=
-g^{\mu\bar\gamma}\na_{\bar k}\na_j \dot g_{\bar\gamma\nu}.
\eea
To apply this formula to the case of the Anomaly flow, where $\p_tg_{\bar\gamma\nu}$ is given by Theorem \ref{anomaly4}, we need to work out the covariant derivatives of the curvature tensor for Hermitian metrics. This is done in the following lemma:

\begin{lemma}
\label{bianchi3}
Let $\o$ be any Hermitian metric (not necessarily conformally balanced). Then we have the following identities
\bea
\na_{\bar k}\na_j R_{\bar\gamma s\bar\mu\lambda}
&=&
\na_s\na_{\bar\gamma}R_{\bar kj\bar\mu\lambda}
+
\na_{\bar k}(T^r{}_{sj}R_{\bar\gamma r\bar\mu\lambda})
+
\na_s(\bar T^{\bar r}{}_{\bar\gamma\bar k}R_{\bar r j\bar\mu\lambda})
\nonumber\\
&&
-R_{\bar ks\bar\gamma}{}^{\bar\kappa}R_{\bar\kappa j\bar\mu\lambda}
+R_{\bar ks}{}^\kappa{}_jR_{\bar\gamma\kappa\bar\mu\lambda}
-R_{\bar ks\bar\mu}{}^{\bar\kappa}R_{\bar\gamma j\bar\kappa\lambda}
+R_{\bar ks}{}^\kappa{}_\lambda R_{\bar\gamma j\bar\mu\kappa}.
\nonumber\\
\na_{\bar k}\na_j\tilde R_{\bar\mu\lambda}
&=&
\Delta R_{\bar kj\bar\mu\lambda}
+
\na_{\bar k}(T^r{}_{s j}R{}^s{}_{r\bar\mu\lambda})
+
\na^{\bar\gamma}(\bar T^{\bar r}{}_{\bar\gamma\bar k}R_{\bar rj\bar\mu\lambda})
\nonumber\\
&&
-R'_{\bar k}{}^{\bar\kappa}R_{\bar\kappa j\bar\mu\lambda}
+R_{\bar k s}{}^\kappa{}_jR^s{}_{\kappa\bar\mu\lambda}
-R_{\bar ks\bar\mu}{}^{\bar\kappa}R^s{}_{j\bar\kappa\lambda}
+R_{\bar ks}{}^\kappa{}_\lambda
R^s{}_{j\bar\mu\kappa})
\nonumber\\
\na_{\bar k}\na_j\tilde R
&=&
\Delta R_{\bar kj}+\na_{\bar k}(T^r{}_{sj}R^s{}_r)
+
\na^{\bar\gamma}(\bar T^{\bar r}{}_{\bar\gamma\bar k}R_{\bar rj})
-R'_{\bar k}{}^{\bar\kappa}R_{\bar\kappa j}+R_{\bar ks}{}^\kappa{}_jR^s{}_\kappa.
\eea
\end{lemma}

\medskip

To clarify the notation: we are writing $\Delta = g^{j\bar k} \nabla_j \na_{\bar k}$ for the `rough' Laplacian and $\bar\Delta = g^{j\bar k} \nabla_{\bar k} \na_{j}$ for its conjugate. While $\Delta$ and $\bar\Delta$ agree when acting on functions, they differ by curvature terms when acting on tensors.

\bigskip

\noindent
{\it Proof.} The proof is a straightforward application of the Bianchi identity, beginning with
\bea
\na_{\bar k}\na_j R_{\bar\gamma s\bar\mu\lambda}
=
\na_{\bar k}(\na_s R_{\bar\gamma j\bar\mu\lambda}
+T^r{}_{sj}R_{\bar\gamma r\bar\mu\lambda})
\eea
and applying it again, after commuting the covariant derivatives $\na_{\bar k}$ and $\na_s$. Q.E.D.

\

We return now to the Anomaly flow of conformally balanced metrics. First, we write
\bea
\p_tR_{\bar kj}{}^\rho{}_\lambda
&=&
-\na_{\bar k}\na_j({1\over 2\|\Omega\|_{\o}}g^{\rho\bar\mu}\Psi_{\bar\mu\lambda})
\nonumber\\
&=&
-{1\over 2\|\Omega\|_{\o}}g^{\rho\bar\mu}\na_{\bar k}\na_j\Psi_{\bar\mu\lambda}
-\na_{\bar k}({1\over 2\|\Omega\|_{\o}})\na_j(g^{\rho\bar\mu}\Psi_{\bar\mu\lambda})
-\na_j({1\over 2\|\Omega\|_{\o}})\na_{\bar k}(g^{\rho\bar\mu}\Psi_{\bar\mu\lambda})
\nonumber\\
&&
-\na_{\bar k}\na_j({1\over 2\|\Omega\|_{\o}})g^{\rho\bar\mu}\Psi_{\bar\mu\lambda}
\nonumber\\
&=&
-{1\over 2\|\Omega\|_{\o}}g^{\rho\bar\mu}\na_{\bar k}\na_j\Psi_{\bar\mu\lambda}
+{1\over 2\|\Omega\|_\o}\bar T_{\bar k}\na_j\Psi^\rho{}_\lambda
+{1\over 2\|\Omega\|_\o}T_{j}\na_{\bar k}\Psi^\rho{}_\lambda
\nonumber\\
&&
+{1\over 2\|\Omega\|_\o}({1\over 2}R_{\bar kj}-T_j\bar T_{\bar k})\Psi^\rho{}_\lambda
\eea
where we used (iii) in Lemma \ref{cb} to get the last equality.

We concentrate on the first term, which can be written in the following way, using Lemma \ref{Psi1},
\bea
-{1\over 2\|\Omega\|_{\o}}g^{\rho\bar\mu}\na_{\bar k}\na_j\Psi_{\bar\mu\lambda}
&=&
{1\over 2\|\Omega\|_{\o}}
g^{\rho\bar\mu}\na_{\bar k}\na_j
\tilde R_{\bar \mu\lambda}+
{1\over 2\|\Omega\|_\o}g^{\rho\bar\mu}g^{s\bar r}\alpha'\na_{\bar k}\na_j
(R_{[\bar\mu s}{}^\alpha{}_\beta R_{\bar r\lambda]}{}^\beta{}_\alpha)
\nonumber\\
&&
\quad
-{1\over 2\|\Omega\|_{\o}}
g^{\rho\bar\mu}\na_{\bar k}\na_j
((T\bar T)_{\bar\mu\lambda}+\alpha'\Phi_{\bar\mu\lambda}).
\eea
The terms in the second line are lower order terms that we shall leave as they are for the moment, and just collect them at the end. 
The first term on the right hand side can be rewritten as follows, using Lemma
\ref{bianchi3},
\bea
{1\over 2\|\Omega\|_\o}
g^{\rho\bar\mu}\na_{\bar k}\na_j\tilde R_{\bar\mu\lambda}
&=&
{1\over 2\|\Omega\|_\o}
\Delta R_{\bar kj}{}^\rho{}_\lambda
+
{1\over 2\|\Omega\|_\o}
\bigg[\na_{\bar k}(T^r{}_{s j}R{}^s{}_{r}{}^\rho
{}_\lambda)
+
\na^{\bar\gamma}(\bar T^{\bar r}{}_{\bar\gamma\bar k}R_{\bar rj}{}^\rho{}_\lambda)
\nonumber\\
&&
-R'_{\bar k}{}^{\bar\kappa}R_{\bar\kappa j}{}^\rho{}_\lambda
+R_{\bar k s}{}^\kappa{}_jR^s{}_{\kappa}{}^\rho{}_\lambda
-R_{\bar ks}{}^{\rho\bar\kappa}R^s{}_{j\bar\kappa\lambda}
+R_{\bar ks}{}^\kappa{}_\lambda
R^s{}_{j}{}^\rho{}_\kappa\bigg].
\eea
It remains only to work out the contribution of the second term on the right hand side,
\bea
&&{1\over 2\|\Omega\|_\o}\alpha' g^{\rho\bar\mu}g^{s\bar r}
\na_{\bar k}\na_j(R_{[\bar\mu s}{}^\alpha{}_\beta R_{\bar r\lambda]}{}^\beta{}_\alpha)\\\nonumber
&=&
{\alpha' g^{\rho\bar\mu}g^{s\bar r}\over 2\|\Omega\|_\o}2(\na_{\bar k}\na_j R_{[\bar\mu s}{}^\alpha{}_\beta)
R_{\bar r\lambda]}{}^\beta{}_\alpha
+
{\alpha'g^{\rho\bar\mu}g^{s\bar r}\over 2\|\Omega\|_\o}
(\na_j R_{[\bar\mu s}{}^\alpha{}_\beta \na_{\bar k}R_{\bar r\lambda]}{}^\beta{}_\alpha
+
\na_{\bar k}R_{[\bar\mu s}{}^\alpha{}_\beta \na_jR_{\bar r\lambda]}{}^\beta{}_\alpha).\nonumber
\eea
Again the second term on the right hand side contains only lower order terms, which we leave as they are and collect only at the end. Using Lemma \ref{bianchi3}, the first term can be rewritten as,
\bea
{\alpha' g^{\rho\bar\mu}g^{s\bar r}\over 2\|\Omega\|_\o}2(\na_{\bar k}\na_j R_{[\bar\mu s}{}^\alpha{}_\beta)
R_{\bar r\lambda]}{}^\beta{}_\alpha
&=&
{\alpha'g^{\rho\bar\mu}g^{s\bar r}\over 2\|\Omega\|_\o}
2
R_{[\bar r\lambda}{}^{\beta\bar\delta}
\na_s\na_{\bar\mu]}R_{\bar kj\bar\delta\beta}
\\
&&
+
{\alpha' g^{\rho\bar\mu}g^{s\bar r}\over 2\|\Omega\|_\o}2R_{[\bar r[\lambda}{}^{\beta\bar\delta}
\bigg[
\na_{\bar k}(T^r{}_{s]j}R_{\bar\mu] r\bar\delta\beta})
+
\na_{s]}(\bar T^{\bar r}{}_{\bar\mu]\bar k}R_{\bar r j\bar\delta\beta})
\nonumber\\
&&
-R_{\bar ks]\bar\mu]}{}^{\bar\kappa}R_{\bar\kappa j\bar\delta\beta}
+R_{\bar ks]}{}^\kappa{}_jR_{\bar\mu]\kappa\bar\delta\beta}
\nonumber\\
&&
-R_{\bar ks]\bar\delta}{}^{\bar\kappa}R_{\bar\mu] j\bar\kappa\beta}
+R_{\bar ks]}{}^\kappa{}_\beta R_{\bar\mu] j\bar\delta\kappa}
\bigg]
\nonumber
\eea
where we have again anti-symmetrized in the unbarred indices $s$ and $\lambda$, and separately in the barred indices $\bar\mu$ and $\bar r$. Whenever there are many indices in the same row and whenever a more explicit designation may be helpful, we have indicated the indices to be anti-symmetrized, either by a symbol $[$ on the left or a symbol $]$ on the right of the relevant index.

\medskip

We obtain in this way the following theorem:

\begin{theorem}
\label{curvatureflow}
Consider the Anomaly flow (\ref{anomaly1}) with an initial metric $\o_0$ which is conformally balanced. Then the curvature of the metric flows according to the following equation
\bea
\p_t R_{\bar kj}{}^\rho{}_\lambda
&=&
{1\over 2\|\Omega\|_\o}\big(\Delta R_{\bar kj}{}^\rho{}_\lambda
+
2\alpha' g^{\rho\bar\mu}g^{s\bar r}R_{[\bar r\lambda}{}^\beta{}_\alpha
\na_s\na_{\bar\mu]} R_{\bar kj}{}^\alpha{}_\beta\big)
\nonumber\\
&&
+{1\over 2\|\Omega\|_\o}\bar T_{\bar k}\na_j\Psi^\rho{}_\lambda
+{1\over 2\|\Omega\|_\o}T_{j}\na_{\bar k}\Psi^\rho{}_\lambda
+{1\over 2\|\Omega\|_\o}({1\over 2}R_{\bar kj}-T_j\bar T_{\bar k})\Psi^\rho{}_\lambda
\nonumber\\
&&
-{1\over 2\|\Omega\|_{\o}}
g^{\rho\bar\mu}\na_{\bar k}\na_j
((T\bar T)_{\bar\mu\lambda}+\alpha'\Phi_{\bar\mu\lambda})
\nonumber\\
&&
+{1\over 2\|\Omega\|_\o}
\bigg[\na_{\bar k}(T^r{}_{s j}R{}^s{}_{r}{}^\rho
{}_\lambda)
+
\na^{\bar\gamma}(\bar T^{\bar r}{}_{\bar\gamma\bar k}R_{\bar rj}{}^\rho{}_\lambda)
\nonumber\\
&&
\qquad\qquad
-R'_{\bar k}{}^{\bar\kappa}R_{\bar\kappa j}{}^\rho{}_\lambda
+R_{\bar k s}{}^\kappa{}_jR^s{}_{\kappa}{}^\rho{}_\lambda
-R_{\bar ks}{}^{\rho\bar\kappa}R^s{}_{j\bar\kappa\lambda}
+R_{\bar ks}{}^\kappa{}_\lambda
R^s{}_{j}{}^\rho{}_\kappa\bigg]
\nonumber\\
&&
+
{\alpha'g^{\rho\bar\mu}g^{s\bar r}\over 2\|\Omega\|_\o}
(\na_j R_{[\bar\mu s}{}^\alpha{}_\beta \na_{\bar k}R_{\bar r\lambda]}{}^\beta{}_\alpha
+
\na_{\bar k}R_{[\bar\mu s}{}^\alpha{}_\beta \na_jR_{\bar r\lambda]}{}^\beta{}_\alpha)
\nonumber\\
&&
+
{\alpha' g^{\rho\bar\mu}g^{s\bar r}\over 2\|\Omega\|_\o}2R_{[\bar r[\lambda}{}^{\beta\bar\delta}
\bigg[
\na_{\bar k}(T^r{}_{s]j}R_{\bar\mu] r\bar\delta\beta})
+
\na_{s]}(\bar T^{\bar r}{}_{\bar\mu]\bar k}R_{\bar r j\bar\delta\beta})
\nonumber\\
&&
-R_{\bar ks]\bar\mu]}{}^{\bar\kappa}R_{\bar\kappa j\bar\delta\beta}
+R_{\bar ks]}{}^\kappa{}_jR_{\bar\mu]\kappa\bar\delta\beta}
-R_{\bar ks]\bar\delta}{}^{\bar\kappa}R_{\bar\mu] j\bar\kappa\beta}
+R_{\bar ks]}{}^\kappa{}_\beta R_{\bar\mu] j\bar\delta\kappa}
\bigg]
\eea
\end{theorem}

\

\subsection{Flow of the Ricci curvature}

The flow of the Riemann curvature tensor implies immediately that of the Ricci curvature,
\bea
\p_tR_{\bar kj}
&=&
{1\over 2\|\Omega\|_\o}\big(\Delta R_{\bar kj}
+
2\alpha' g^{\lambda\bar\mu}g^{s\bar r}R_{[\bar r\lambda}{}^\beta{}_\alpha
\na_s\na_{\bar\mu]} R_{\bar kj}{}^\alpha{}_\beta
\nonumber\\
&&
+{1\over 2\|\Omega\|_\o}\bar T_{\bar k}\na_j\Psi^\lambda{}_\lambda
+{1\over 2\|\Omega\|_\o}T_{j}\na_{\bar k}\Psi^\lambda{}_\lambda
+{1\over 2\|\Omega\|_\o}({1\over 2}R_{\bar kj}-T_j\bar T_{\bar k})\Psi^\lambda{}_\lambda
\nonumber\\
&&
-{1\over 2\|\Omega\|_{\o}}
\na_{\bar k}\na_j
(|T|^2+\alpha'\Phi^\lambda{}_{\lambda})
\nonumber\\
&&
+{1\over 2\|\Omega\|_\o}
\bigg[\na_{\bar k}(T^r{}_{s j}R{}^s{}_{r})
+
\na^{\bar\gamma}(\bar T^{\bar r}{}_{\bar\gamma\bar k}R_{\bar rj})
-R'_{\bar k}{}^{\bar\kappa}R_{\bar\kappa j}{}
+R_{\bar k s}{}^\kappa{}_jR^s{}_{\kappa}\bigg]
\nonumber\\
&&
+
{\alpha'g^{\lambda\bar\mu}g^{s\bar r}\over 2\|\Omega\|_\o}
(\na_j R_{[\bar\mu s}{}^\alpha{}_\beta \na_{\bar k}R_{\bar r\lambda]}{}^\beta{}_\alpha
+
\na_{\bar k}R_{[\bar\mu s}{}^\alpha{}_\beta \na_jR_{\bar r\lambda]}{}^\beta{}_\alpha)
\nonumber\\
&&
+
{\alpha' g^{\lambda\bar\mu}g^{s\bar r}\over 2\|\Omega\|_\o}2R_{[\bar r[\lambda}{}^{\beta\bar\delta}
\bigg[
\na_{\bar k}(T^r{}_{s]j}R_{\bar\mu] r\bar\delta\beta})
+
\na_{s]}(\bar T^{\bar r}{}_{\bar\mu]\bar k}R_{\bar r j\bar\delta\beta})
\nonumber\\
&&
-R_{\bar ks]\bar\mu]}{}^{\bar\kappa}R_{\bar\kappa j\bar\delta\beta}
+R_{\bar ks]}{}^\kappa{}_jR_{\bar\mu]\kappa\bar\delta\beta}
-R_{\bar ks]\bar\delta}{}^{\bar\kappa}R_{\bar\mu] j\bar\kappa\beta}
+R_{\bar ks]}{}^\kappa{}_\beta R_{\bar\mu] j\bar\delta\kappa}
\bigg]
\eea
with $|T|^2= g^{j\bar k}g^{s\bar r} g^{m\bar \ell} T_{\bar r m j} \bar T_{s\bar \ell \bar k}$.

\

\subsection{Flow of the scalar curvature}

If we write $R=g^{j\bar k}R_{\bar kj}$, we obtain
\bea
\p_tR=g^{j\bar k}\p_tR_{\bar kj}-g^{j\bar m}\p_t g_{\bar mq}g^{q\bar k}R_{\bar kj}.
\eea
Applying the preceding formula for the flow $\p_tR_{\bar kj}$ of the Ricci curvature, we find
\bea
\p_tR
&=&
{1\over 2\|\Omega\|_\o}\big(\Delta R
+
2\alpha' g^{\lambda\bar\mu}g^{s\bar r}R_{[\bar r\lambda}{}^\beta{}_\alpha
\na_s\na_{\bar\mu]} \tilde R^\alpha{}_\beta\big)
\nonumber\\
&&
+{1\over 2\|\Omega\|_\o}
\big(\bar T^j\na_j\Psi^\lambda{}_\lambda
+{1\over 2\|\Omega\|_\o}T^{\bar k}\na_{\bar k}\Psi^\lambda{}_\lambda
+({1\over 2}R-T_j \bar{T}^j)\Psi^\lambda{}_\lambda\big)
\nonumber\\
&&
-{1\over 2\|\Omega\|_{\o}}
\Delta
(|T|^2+\alpha'\Phi^\lambda{}_{\lambda})
-
{1\over 2\|\Omega\|_\o}R^{q\bar m}\Psi_{\bar mq}
\nonumber\\
&&
+{1\over 2\|\Omega\|_\o}
\bigg(\na_{\bar k}(T^r{}_{s j}R{}^s{}_{r})
+
\na^{\bar\gamma}(\bar T^{\bar r}{}_{\bar\gamma\bar k}R_{\bar rj})\bigg)
\nonumber\\
&&
+
{\alpha'g^{\lambda\bar\mu}g^{s\bar r}\over 2\|\Omega\|_\o}
(\na_j R_{[\bar\mu s}{}^\alpha{}_\beta \na^jR_{\bar r\lambda]}{}^\beta{}_\alpha
+
\na_{\bar k}R_{[\bar\mu s}{}^\alpha{}_\beta \na^{\bar k}R_{\bar r\lambda]}{}^\beta{}_\alpha)
\nonumber\\
&&
+
{\alpha' g^{\lambda\bar\mu}g^{s\bar r}\over 2\|\Omega\|_\o}2R_{[\bar r[\lambda}{}^{\beta\bar\delta}
\bigg[
\na^j(T^\gamma{}_{s]j}R_{\bar\mu] \gamma\bar\delta\beta})
+
\na_{s]}(\bar T^{\bar \gamma}{}_{\bar\mu]}{}^jR_{\bar \gamma j\bar\delta\beta})
\nonumber\\
&&
-R^j{}_{s]\bar\mu]}{}^{\bar\kappa}R_{\bar\kappa j\bar\delta\beta}
+\tilde R_{s]}{}^\kappa R_{\bar\mu]\kappa\bar\delta\beta}
-R^j{}_{s]\bar\delta}{}^{\bar\kappa}R_{\bar\mu] j\bar\kappa\beta}
+R^j{}_{s]}{}^\kappa{}_\beta R_{\bar\mu] j\bar\delta\kappa}
\bigg].
\eea

\

\subsection{Flow of the torsion tensor}

We differentiate the coefficients $T_{\bar pjq}$ of the torsion tensor,
\bea
\p_tT_{\bar pjq}&=&\p_j\dot g_{\bar pq}-\p_q\dot g_{\bar pj}
\nonumber\\
&=&
\p_j({1\over 2\|\Omega\|_\o}\Psi_{\bar pq})-\p_q({1\over 2\|\Omega\|_\o}\Psi_{\bar pj})
\nonumber\\
&=&
{1\over 2\|\Omega\|_\o}(\na_j\Psi_{\bar pq}-\na_q\Psi_{\bar pj}+T^m{}_{jq}\Psi_{\bar pm})
-{1\over 2\|\Omega\|_\o}(T_j\Psi_{\bar pq}-T_q\Psi_{\bar pj}).
\eea
Once again, we concentrate on the leading term, which is
\bea
{1\over 2\|\Omega\|_\o}
(\na_j\Psi_{\bar pq}-\na_q\Psi_{\bar pj})
&=&
{1\over 2\|\Omega\|_\o}(\na_j(-\tilde R_{\bar pq}+(T\bar T)_{\bar pq})
-\na_q(-\tilde R_{\bar pj}+(T\bar T)_{\bar pj})
\nonumber\\
&&
-
{1\over 2\|\Omega\|_\o}\alpha'
g^{s\bar r}\na_j (R_{[\bar p s}{}^\alpha{}_\beta R_{\bar rq]}{}^\beta{}_\alpha
-\Phi_{\bar p s\bar rq})
\nonumber\\
&&
+{1\over 2\|\Omega\|_\o}\alpha'
g^{s\bar r} \na_q(R_{[\bar p s}{}^\alpha{}_\beta R_{\bar rj]}{}^\beta{}_\alpha)
-
\Phi_{\bar ps\bar rj}).
\eea
Although this is not apparent at first sight, the key diffusion term $\Delta T_{\bar pjq}$ can be extracted from the right hand side. 
The basic identity in this case is the following:

\begin{lemma}
Let $\o$ be any Hermitian metric (not necessarily conformally balanced). Then
\bea
(\Delta T)_{\bar p jq}
=
\na_q\tilde R_{\bar pj}-\na_j\tilde R_{\bar pq}
+
T^r{}_{q\lambda}R^\lambda{}_{r\bar pj}
-
T^r{}_{j\lambda}R^\lambda{}_{r\bar p q}.
\eea
\end{lemma}

\medskip
\noindent
{\it Proof.} We compute the components of the left hand side, using the Bianchi identities,
\bea
(\Delta T)_{\bar pjq}
&=&
g^{\lambda\bar\mu}\na_\lambda\na_{\bar\mu}T_{\bar pjq}
\nonumber\\
&=&
g^{\lambda\bar\mu}\na_{\lambda}(R_{\bar \mu q\bar p j}
-
R_{\bar\mu j\bar pq})
\nonumber\\
&=&
g^{\lambda\bar\mu}(\na_q R_{\bar\mu\lambda\bar pj}-\na_jR_{\bar\mu\lambda\bar pq}+T^r{}_{q\lambda}R_{\bar\mu r\bar pj}
-
T^r{}_{j\lambda}R_{\bar\mu r\bar p q}).
\nonumber\\
&=&
\na_q\tilde R_{\bar pj}-\na_j\tilde R_{\bar pq}+T^r{}_{q\lambda}R^\lambda{}_{r\bar pj}
-
T^r{}_{j\lambda}R^\lambda{}_{r\bar p q}.
\eea
This proves the lemma.

\bigskip
Comparing this identity with the previous expression that we derived for $\p_tT_{\bar pjq}$, we obtain the following theorem:

\begin{theorem}
\label{torsionflow}
Consider the Anomaly flow (\ref{anomaly1}) with an initial metric $\o_0$ which is conformally balanced. Then the flow of the torsion $T=i\p\o$ is given by
\bea
\p_tT_{\bar pjq}&=&{1\over 2\|\Omega\|_\o}
\bigg[\Delta T_{\bar pjq}
-\alpha'
g^{s\bar r}(\na_j (R_{[\bar p s}{}^\alpha{}_\beta R_{\bar rq]}{}^\beta{}_\alpha
-\Phi_{\bar p s\bar rq})
+\alpha'
g^{s \bar{r}} \na_q(R_{[\bar p s}{}^\alpha{}_\beta R_{\bar rj]}{}^\beta{}_\alpha
-
\Phi_{\bar ps\bar rj}))\bigg]
\nonumber\\
&&
+{1\over 2\|\Omega\|_\o}(T^m{}_{jq}\Psi_{\bar pm}-T_j\Psi_{\bar pq}+T_q\Psi_{\bar pj}+\na_j (T\bar T)_{\bar pq}-\na_q(T\bar T)_{\bar pj})
\nonumber\\
&&
-{1\over 2\|\Omega\|_\o}
(T^r{}_{q\lambda}R^\lambda{}_{r\bar pj}
-
T^r{}_{j\lambda}R^\lambda{}_{r\bar p q})
\eea
\end{theorem}

\

\

\section{A Model Problem: $\alpha'=0$}
\setcounter{equation}{0}

A first model which is simpler than the full Anomaly flow and whose study could be instructive, is obtained by setting $\alpha'=0$. While this special case eliminates the quadratic terms in the curvature tensor in (\ref{anomaly1}), it still presents some new difficulties relative to the well-known Ricci flow and Donaldson heat flow because of the evolving torsion. More precisely, we shall consider the flow
\bea
\label{anomalyK}
\p_t(\|\Omega\|_\o\o^2)=i\p\bar\p\o.
\eea
The stationary points of the flow are then given by the equivalent equations
\bea
i\p\bar\p\o=0
\eea
for a Hermitian metric satisfying the conformally balanced condition $d(\|\Omega\|_{\o}\o^2)=0$. Such a Hermitian metric must be K\"ahler and Ricci-flat \cite{IP}, since contracting (\ref{iddb-omega-identity}) and applying Lemma \ref{torsion/curvature} shows that such a metric must satisfy
\be
g^{j \bar{k}} \p_j \p_{\bar{k}} \log \| \Omega \|^2 = \tilde{R} = |T|^2.
\ee
By the maximum principle, $|T|^2=0$ and $\log \| \Omega \|^2$ is constant. Thus the Anomaly flow with $\alpha'=0$ can be used to determine whether a conformally balanced manifold is actually K\"ahler. 

We note that the flow (\ref{anomalyK}) is also related to the problem of prescribing metrics in a balanced class raised in the recent survey of Garcia-Fernandez \cite{Ga}, so our results in this section can be viewed as a first step towards an eventual solution.

\medskip

\subsection{Flow of the curvature and the torsion}

For convenience, we summarize here the main formulas for the Anomaly flow (\ref{anomalyK}). They can be obtained from the general formulas obtained earlier by setting $\alpha'=0$. Let us still use $\Psi$ to denote the right hand side of the flow, that is $\Psi_{\bar p q} = -\tilde R_{\bar pq} + (T\bar T)_{\bar pq}$.
Then the flow of the metric is given by
\bea
\p_t g_{\bar pq}
= {1\over 2\|\Omega\|_\o} \Psi_{\bar p q} =
{1\over 2\|\Omega\|_\o}\bigg[-\tilde R_{\bar pq}
+
(T\bar T)_{\bar pq}\bigg]
\eea
while the flows of the curvature tensors are given by
\bea
\p_t R_{\bar kj}{}^\rho{}_\lambda
&=&
{1\over 2\|\Omega\|_\o}\Delta R_{\bar kj}{}^\rho{}_\lambda
-
{1\over 2\|\Omega\|_{\o}}
g^{\rho\bar\mu}\na_{\bar k}\na_j
(T\bar T)_{\bar\mu\lambda} 
\nonumber\\
&&
+{1\over 2\|\Omega\|_\o}\bigg(\bar T_{\bar k}\na_j
+T_{j}\na_{\bar k}
+({1\over 2}R_{\bar kj}-T_j\bar T_{\bar k})\bigg)\Psi^\rho{}_\lambda
\nonumber\\
&&
+{1\over 2\|\Omega\|_\o}
\bigg[\na_{\bar k}(T^r{}_{s j}R{}^s{}_{r}{}^\rho
{}_\lambda)
+
\na^{\bar\gamma}(\bar T^{\bar r}{}_{\bar\gamma\bar k}R_{\bar rj}{}^\rho{}_\lambda)
\nonumber\\
&&
-R'_{\bar k}{}^{\bar\kappa}R_{\bar\kappa j}{}^\rho{}_\lambda
+R_{\bar k s}{}^\kappa{}_jR^s{}_{\kappa}{}^\rho{}_\lambda
-R_{\bar ks}{}^{\rho\bar\kappa}R^s{}_{j\bar\kappa\lambda}
+R_{\bar ks}{}^\kappa{}_\lambda
R^s{}_{j}{}^\rho{}_\kappa\bigg]
\nonumber\\
&&
\nonumber\\
\p_tR_{\bar kj}
&=&
{1\over 2\|\Omega\|_\o}\Delta R_{\bar kj}-
{1\over 2\|\Omega\|_{\o}}
\na_{\bar k}\na_j
|T|^2
\nonumber\\
&&
+{1\over 2\|\Omega\|_\o}\bigg(\bar T_{\bar k}\na_j
+T_{j}\na_{\bar k}
+({1\over 2}R_{\bar kj}-T_j\bar T_{\bar k})\bigg)
(-R+|T|^2)
\nonumber\\
&&
+{1\over 2\|\Omega\|_\o}
\bigg[\na_{\bar k}(T^r{}_{s j}R{}^s{}_{r})
+
\na^{\bar\gamma}(\bar T^{\bar r}{}_{\bar\gamma\bar k}R_{\bar rj})
-R'_{\bar k}{}^{\bar\kappa}R_{\bar\kappa j}{}
+R_{\bar k s}{}^\kappa{}_jR^s{}_{\kappa}\bigg]
\nonumber\\
\p_tR
&=&
{1\over 2\|\Omega\|_\o}\Delta R-
{1\over 2\|\Omega\|_{\o}}
\Delta
|T|^2-{1\over 2\|\Omega\|_\o}R^{j\bar k}\Psi_{\bar k j}
\nonumber\\
&&
+{1\over 2\|\Omega\|_\o}\bigg(\bar T_{\bar k}\na^{\bar k}
+T_{j}\na^j
+({1\over 2}R-T_j \bar{T}^j)\bigg)
(-R+|T|^2)
\nonumber\\
&&
+{1\over 2\|\Omega\|_\o}
\bigg(\na_{\bar k}(T^r{}_{s j}R{}^s{}_{r})
+
\na^{\bar\gamma}(\bar T^{\bar r}{}_{\bar\gamma\bar k}R_{\bar rj})\bigg)
\eea
and the flow of the torsion is given by
\bea
\p_tT_{\bar pjq}
&=&
{1\over 2\|\Omega\|_\o}
\Delta T_{\bar pjq}
-{1\over 2\|\Omega\|_\o}
(T^r{}_{q\lambda}R^\lambda{}_{r\bar pj}
-
T^r{}_{j\lambda}R^\lambda{}_{r\bar p q})
\nonumber\\
&&
+{1\over 2\|\Omega\|_\o}(T^m{}_{jq}\Psi_{\bar pm}-T_j\Psi_{\bar pq}+T_q\Psi_{\bar pj}+\na_j (T\bar T)_{\bar pq}-\na_q(T\bar T)_{\bar pj}).
\eea
For later use, we also record here the flow of the norm $\|\Omega\|_\o$,
\bea
\p_t\|\Omega\|_\o
=
{1\over 4}(R-|T|^2 ).
\eea

\

\subsection{Estimates for derivatives of curvature and torsion}

The goal in this section is to prove Theorem \ref{Higherorderestimate}.
We shall use $D$ to denote the derivative when we do not distinguish between $\nabla$ and $\bar\nabla$.
For example, $|D T|$ would include both $|\nabla T|$ and $|\bar\nabla T|$,
and
\bea
|D^kT|^2=\sum_{i+j =k}|\nabla^i\bar\nabla^{j}T|^2.
\eea

\medskip

The proof of Theorem \ref{Higherorderestimate} is by induction on $k$.  The idea is find a suitable test function $G_k(z, t)$ for each $k$, similar to the Ricci flow, and apply the maximum principle.

\medskip

We will first prove the estimate (\ref{CTestimate}) for $k=1$ case. Then, we assume that, for any $0\leq j \leq k-1$, 
\bea\label{assumption1}
|D^j Rm(z, t)|_{\o} \leq {C_j A\over t^{j/2}}, \ \ \ |D^{j+1} T(z, t)|_{\o} \leq  {C_j A\over t^{j/2}}
\eea
for all $z\in M$ and $t\in (0, {1\over A}]$ and show the estimate also holds for $j=k$.

\

We already have the flows of the curvature and of the torsion, as given above. To prove the theorem, we shall also need the flows of their covariant derivatives. They are given in the following lemmas.

\begin{lemma}
\label{DkRm}
Under the induction assumption (\ref{assumption1}) and $|T|^2 \leq A$, we have
\bea
\p_t |D^k Rm|^2
&\leq & 
{1\over 2\|\Omega\|_{\o}}\, 
\Bigg\{ {1\over 2}\Delta_{\R} |D^k Rm|^2  - {3\over 4}|D^{k+1} Rm|^2 
\\
&&
+ CA^{1\over 2} \left(|D^{k+1} Rm|+ |D^{k+2} T| \right) \cdot |D^k Rm|
 \nonumber\\
&&
+ CA \left(|D^{k} Rm|
+ |D^{k+1} T| \right) \cdot |D^k Rm| \nonumber\\
&&+C\, A^2\, t^{-{k\over 2}} \cdot |D^k Rm| + C A^3 t^{-k}
\Big\}\nonumber
\eea
where we write $\Delta_{\R}=\Delta + \bar\Delta$ and $\Delta = g^{\bar q p}\na_{p}\na_{\bar q}$.
\end{lemma}

\medskip
\noindent
{\it Proof.} First, we observe that the flow of the curvature tensor can be expressed as
\bea
 \p_t Rm 
 &=&
 {1\over 2\|\Omega\|_{\o}}\,
 \Big\{
 {1\over 2}\Delta_{\R} Rm + \nabla\bar\nabla(T\ast \bar T) + \bar\nabla(T\ast Rm) + \nabla(\bar T \ast Rm)\\
  \nonumber
  &&
  + Rm\ast Rm + (\bar\nabla T - \bar T \ast T )\ast \Psi + \bar T \ast \nabla \Psi + T\ast \bar\nabla \Psi 
 \Big\}.
\eea
To clarify notation: if $E$ and $F$ are tensors, we write $E\ast F$ for any linear combination of products of the tensors $E$ and $F$ formed by contractions on $E_{i_1\cdots i_k}$ and $F_{j_1\cdots j_l}$ using the metric $g$.

\medskip

Let the terms in the large bracket be denoted by $H$, that is
 \bea
 \p_t Rm = {1\over 2\|\Omega\|_{\o}} H.
 \eea

In general, the Chern unitary connection of a Hermitian metric $g_{\bar kj}$ evolves by
\bea
\p_t A_{\bar km}^j=0,
\qquad
\p_t A_{km}^j=g^{j\bar p}\na_k(\p_t g_{\bar pm}).
\eea
This implies
\bea
\p_t (\nabla^m \bar \nabla^{\ell} Rm) 
= \nabla^m \bar \nabla^{\ell} (\p_t Rm) + \sum_{i+j>0}\sum_{i=0}^m \sum_{j=0}^\ell \nabla^{m-i}\bar\nabla^{\ell - j} Rm \ast \nabla^i \bar \nabla^j (\p_t g)
\eea
Using the evolution equation of $Rm$, we get
\bea
\p_t (\nabla^m \bar \nabla^{\ell} Rm) 
&=&
 \sum_{i=1}^m \sum_{j=1}^\ell \nabla^{m-i}\bar\nabla^{\ell - j} Rm \ast \nabla^i \bar \nabla^j (\p_t g)\\\nonumber
 &&+{1\over 2\|\Omega\|_{\o}}\,\nabla^m \bar\nabla^\ell H + \sum_{i+j>0}\sum_{i=0}^m \sum_{j=0}^\ell \nabla^{m-i} \bar\nabla^{\ell-j} H \ast \nabla^i \bar \nabla^j \left( {1\over 2\|\Omega\|_{\o}}\right)
\eea
We compute the second term,
\bea
\nabla^m \bar\nabla^\ell H 
&=& 
{1\over 2} \nabla^m \bar\nabla^\ell \Delta_{\R} Rm + \nabla^m \bar\nabla^\ell \nabla\bar\nabla (T\ast T) + \nabla^m \bar \nabla^{\ell+1} (T\ast Rm) 
\nonumber\\
&&+ \nabla^m \bar \nabla^\ell \nabla (\bar T \ast Rm) + \nabla^m \bar \nabla^\ell (Rm \ast Rm) + \nabla^m \bar \nabla^{\ell +1} (T\ast \Psi) 
\nonumber\\
&&+ \nabla^m \bar \nabla^\ell (\Psi \ast \bar T\ast T) + \nabla^m \bar \nabla^\ell ( \nabla\Psi \ast \bar T)+ \nabla^m \bar \nabla^\ell ( T \ast \bar{\na} \Psi) 
\eea
In view of the commutation identity given in the appendix,
\bea
\nabla^m\bar \nabla^\ell \Delta_{\R} Rm 
&=&
\Delta_{\R} (\nabla^m \bar \nabla^\ell Rm) + \sum_{i=0}^m \sum_{j=0}^\ell \nabla^i \bar \nabla^j Rm \ast \nabla^{m-i} \bar \nabla^{\ell - j} Rm
\\\nonumber
&&+ \sum_{i=0}^m \sum_{j=0}^\ell \nabla^i \bar \nabla^j T \ast \nabla^{m-i} \bar \nabla^{\ell +1 -i} Rm + \sum_{i=0}^m \sum_{j=0}^\ell \nabla^i \bar \nabla^j T \ast \nabla^{m+1-i} \bar \nabla^{\ell -j} Rm
\eea
we obtain
\bea
\p_t (\nabla^m \bar \nabla^{\ell} Rm) 
&=&
{1\over 2\|\Omega\|_{\o}}\,
 \Big\{
{1\over 2}\Delta_{\R}(\nabla^m \bar \nabla^{\ell} Rm)+ \sum_{i=0}^m \sum_{j=0}^\ell \nabla^i \bar \nabla^j Rm \ast \nabla^{m-i} \bar \nabla^{\ell - j} Rm
\nonumber\\
&&+ \sum_{i=0}^m \sum_{j=0}^\ell \nabla^i \bar \nabla^j T \ast \left(\nabla^{m-i} \bar \nabla^{\ell +1 -i} Rm + \nabla^{m+1-i} \bar \nabla^{\ell -j} Rm\right)
+ \nabla^m \bar\nabla^\ell \nabla\bar\nabla (T\ast T) 
\nonumber\\
&&
+ \nabla^m \bar \nabla^{\ell+1} (T\ast Rm) 
+ \nabla^m \bar \nabla^\ell \nabla (\bar T \ast Rm)
 + \nabla^m \bar \nabla^\ell (Rm \ast Rm) 
\nonumber\\
&&
+ \nabla^m \bar \nabla^{\ell +1} (T\ast \Psi) 
+ \nabla^m \bar \nabla^\ell (\Psi \ast \bar T\ast T) + \nabla^m \bar \nabla^\ell ( \nabla\Psi \ast \bar T)
 \Big\}
  \nonumber\\
 &&
 + \sum_{i+j>0}\sum_{i=0}^m \sum_{j=0}^\ell \nabla^{m-i}\bar\nabla^{\ell - j} Rm \ast \nabla^i \bar \nabla^j (\p_t g)
\nonumber\\
 &&+ \sum_{i+j>0}\sum_{i=0}^m \sum_{j=0}^\ell \nabla^{m-i} \bar\nabla^{\ell-j} H \ast \nabla^i \bar \nabla^j \left( {1\over 2\|\Omega\|_{\o}}\right)
\eea
Next we compute
\bea
\p_t |\nabla^m \bar\nabla^\ell Rm|^2 
&\leq & \< \p_t \nabla^m \bar \nabla^\ell Rm,  \, \nabla^m \bar\nabla^\ell Rm\>
+ \< \nabla^m \bar \nabla^{\ell} Rm, \, \p_t \nabla^m \bar \nabla^\ell Rm\>
\\\nonumber
&& + {C \over 2\|\Omega\|_{\o}} |\nabla^m \bar \nabla^\ell Rm|^2 \cdot |\Psi|
\eea
We also compute
\bea
\Delta_{\R} |\nabla^m \bar \nabla^\ell Rm|^2 
&=& \< \Delta_{\R} \nabla^m \bar \nabla^\ell Rm, \, \nabla^m \bar \nabla^\ell Rm\> + \< \nabla^m \bar \nabla^\ell Rm, \, \Delta_{\R} \nabla^m \bar \nabla^\ell Rm\>
\\\nonumber
&& + 2 |\nabla^{m+1} \bar \nabla^\ell Rm|^2+ 2|\bar \nabla \nabla^m \bar \nabla^\ell Rm|^2
\\\nonumber
&=& 
\< \Delta_{\R} \nabla^m \bar \nabla^\ell Rm, \, \nabla^m \bar \nabla^\ell Rm\> + \< \nabla^m \bar \nabla^\ell Rm, \, \Delta_{\R} \nabla^m \bar \nabla^\ell Rm\>
\\\nonumber
&& + 2 |\nabla^{m+1} \bar \nabla^\ell Rm|^2+ 2 |\nabla^m \bar \nabla^{\ell +1} Rm|^2 
\\\nonumber
&&+ 2 \left(|\bar \nabla \nabla^m \bar \nabla^\ell Rm|^2-  |\nabla^m \bar \nabla^{\ell +1} Rm|^2\right)
\eea
We can estimate the last term by a commutation identity.
\bea
\bar \nabla \nabla^m \bar \nabla^\ell Rm - \nabla^m \bar\nabla\bar \nabla^\ell Rm &=& 
 \sum_{i=0}^{m-1} \nabla^i Rm \ast \nabla^{m-1-i} \bar\nabla^\ell Rm
\eea
It follows that
\bea
&&|\bar \nabla \nabla^m \bar \nabla^\ell Rm|^2-  |\nabla^m \bar \nabla^{\ell +1} Rm|^2\\
&\geq & 
- C| \nabla^m \bar \nabla^{\ell +1} Rm| \cdot \sum_{i=0}^{m-1} |\nabla^i Rm \ast \nabla^{m-1-i}\bar\nabla^\ell Rm|
-C\sum_{i=0}^{m-1} |\nabla^i Rm \ast \nabla^{m-1-i} \bar\nabla^\ell Rm|^2
\nonumber
\\
&\geq & 
- C_1| \nabla^m \bar \nabla^{\ell +1} Rm| \cdot \sum_{i=0}^{m-1} |D^i Rm|\cdot |D^{m+\ell -1-i} Rm| - C\sum_{i=0}^{m-1} |D^i Rm|^2\cdot |D^{m+\ell-1-i} Rm|^2\nonumber
\eea
Putting all the computation together, we arrive at
\bea
&&\p_t |\nabla^m \bar \nabla^\ell Rm|^2 \nonumber\\
&\leq & 
{1\over 2\|\Omega\|_{\o}}\, 
\Bigg\{ 
{1\over 2}\Delta_{\R} |\nabla^m \bar \nabla^\ell Rm|^2  - |\nabla^{m+1}\bar \nabla^\ell Rm|^2 - |\nabla^m \bar \nabla^{\ell+1} Rm|^2
\nonumber\\
&& 
+C_1 | \nabla^m \bar \nabla^{\ell +1} Rm| \cdot \sum_{i=0}^{m-1} |D^i Rm|\cdot |D^{m+\ell -1-i} Rm|
+ C\sum_{i=0}^{m-1} |D^i Rm|^2\cdot |D^{m+\ell-1-i} Rm|^2\nonumber\\
&&
+ C |\nabla^m \bar \nabla^\ell Rm|\cdot 
\Bigg[ \sum_{i=0}^m \sum_{j=0}^\ell |\nabla^i \bar \nabla^j Rm|\cdot |\nabla^{m-i}\bar \nabla^{\ell -j}Rm|
+ |\nabla^{m+1}\bar \nabla^{\ell +1} (T\ast T)|
\nonumber\\
&&
 + \sum_{i=0}^m \sum_{j=0}^{\ell -1} |\nabla^i \bar \nabla^j Rm| \cdot |\nabla^{m-i} \bar\nabla^{\ell -j} (T\ast T)|
+ |\nabla^m \bar \nabla^{\ell +1} (T\ast Rm)|
\nonumber\\
&&
+  \sum_{i=0}^m \sum_{j=0}^\ell |\nabla^i \bar \nabla^j T|\cdot \left(\nabla^{m-i}\bar \nabla^{\ell +1 -j} Rm| + |\nabla^{m+1-i } \bar \nabla^{\ell -j} Rm|\right)
\nonumber\\
&&
+ |\nabla^{m+1} \bar\nabla^{\ell} (\bar T \ast Rm)| + \sum_{i=0}^m \sum_{j=0}^{\ell-1} |\nabla^i\bar \nabla^j Rm|\cdot |\nabla^{m-i} \bar\nabla^{\ell-1-j}(\bar T \ast Rm)|
\nonumber\\
&&
+ |\nabla^m \bar\nabla^\ell (Rm\ast Rm)| + |\nabla^m \bar\nabla^{\ell+1} (T\ast \Psi)|
+ |\nabla^m \bar \nabla^\ell (\Psi \ast \bar T\ast T)|  
+ |\nabla^m \bar \nabla^\ell (\nabla\Psi \ast \bar T)|
\nonumber\\
&&
+ \sum_{i+j>0}\sum_{i=0}^m \sum_{j=0}^\ell |\nabla^{m-i}\bar \nabla^{\ell -j} H|\cdot |\nabla^i \bar \nabla^j \left({1\over 2\|\Omega\|_{\o}}\right) |
\nonumber\\
&&
+ \sum_{i+j>0}\sum_{i=0}^m \sum_{j=0}^\ell |\nabla^{m-i}\bar\nabla^{\ell-j} Rm| \cdot |\nabla^i \bar \nabla^j (\p_t g)|
\Bigg]
\Bigg\}
\nonumber\\
&&
+ {C \over 2\|\Omega\|_{\o}} |\nabla^m \bar \nabla^\ell Rm|^2 \cdot |\Psi|
\eea
where we used commutating identities for terms
$\nabla^m \bar\nabla^\ell \nabla\bar\nabla (T\ast T)$ and $\nabla^m \bar \nabla^\ell \nabla (\bar T \ast Rm)$ in the evolution equation $\p_t \nabla^k \bar\nabla^\ell Rm$.
Next, we use the non-standard notation $D$ introduced at the beginning of this section. Note that, for a tensor $E$,
\bea
 |\nabla^i \bar\nabla^j E| \leq |D^{i+j} E|.
\eea
Let $k=m+\ell$. We have
\bea\nonumber
&&\p_t |\nabla^m \bar \nabla^\ell Rm|^2 \\
&\leq & 
{1\over 2\|\Omega\|_{\o}}\, 
\Bigg\{ 
{1\over 2} \Delta_{\R} |\nabla^m \bar \nabla^\ell Rm|^2  - |\nabla^{m+1}\bar \nabla^\ell Rm|^2 - |\nabla^m \bar \nabla^{\ell+1} Rm|^2
\nonumber\\
&& 
+C_1 | \nabla^m \bar \nabla^{\ell +1} Rm| \cdot \sum_{i=0}^{k-1} |D^i Rm|\cdot |D^{k -1-i} Rm|
+ C\sum_{i=0}^{k-1} |D^i Rm|^2\cdot |D^{k-1-i} Rm|^2\nonumber\\
&&
+ C |\nabla^m \bar \nabla^\ell Rm|\cdot 
\Bigg[ 
\sum_{i=0}^k |D^i Rm|\cdot |D^{k-i}Rm| + \sum_{i=0}^k |D^i T|\cdot |D^{k+1 -i}Rm|
\nonumber\\
&&
+ |D^{k+2}(T\ast T)|  + \sum_{i=0}^{k-1} |D^i Rm|\cdot |D^{k-i}(T\ast T)|
\nonumber\\
&&+ |D^{k+1}(T\ast Rm)|
+ |D^{k+1} (\bar T\ast Rm)| 
+ \sum_{i=0}^{k-1} |D^i Rm|\cdot |D^{k-1-i} (\bar T\ast Rm)| 
\nonumber\\
&&
+|D^k (Rm\ast Rm)| + |D^{k+1} (T\ast \Psi)| + |D^{k}(\Psi\ast T\ast T) | + |D^k(\nabla\Psi\ast T)|\nonumber\\
&&
+ \sum_{i=1}^k |D^{k-i} H |\cdot |D^i \left({1\over 2\|\Omega\|_{\o}}\right)|
+ \sum_{i=1}^k |D^{k-i} Rm|\cdot |D^i (\p_t g)|
\Bigg]
\Bigg\}
\nonumber\\
&&
+ {C \over 2\|\Omega\|_{\o}} |\nabla^m \bar \nabla^\ell Rm|^2 \cdot |\Psi|
\eea
Recall that
\bea
|D^k Rm|^2 &=& \sum_{m+\ell =k} |\nabla^m \bar \nabla^\ell Rm|^2\\
|\nabla^m \bar \nabla^{\ell +1} Rm |&\leq& |D^{k +1} Rm|, \ \  \  |\nabla^m \bar \nabla^{\ell} Rm |\leq |D^{k} Rm|
\eea
and we also have
\bea\nonumber
|D^{k+1} Rm|^2 &=& \sum_{m+q =k+1}|\nabla^{m} \bar \nabla^q Rm|^2
=
\sum_{m+q - 1 =k,\, q\geq 1} |\nabla^{m} \bar \nabla^q Rm|^2 + |\nabla^{k+1} Rm|^2
\\\nonumber
&=& \sum_{m+\ell =k, \, m\geq 0, \, \ell \geq 0}  |\nabla^{m} \bar \nabla^{\ell +1} Rm|^2 + |\nabla^{k+1} Rm|^2
\\
&\leq &
\sum_{m+\ell =k}  |\nabla^{m} \bar \nabla^{\ell +1} Rm|^2 +\sum_{m+\ell =k}  |\nabla^{m+1} \bar \nabla^{\ell } Rm|^2 
\eea
Using these inequalities, we get
\bea\label{evolutionDkRm}
\nonumber
&&\p_t |D^k Rm|^2\\\nonumber
&\leq & 
{1\over 2\|\Omega\|_{\o}}\, 
\Bigg\{ 
{1\over 2} \Delta_{\R} |D^k Rm|^2  - |D^{k+1} Rm|^2 \nonumber\\
&& 
+C_1 | D^{k +1} Rm| \cdot \sum_{i=0}^{k-1} |D^i Rm|\cdot |D^{k -1-i} Rm| + C\sum_{i=0}^{k-1} |D^i Rm|^2\cdot |D^{k-1-i} Rm|^2\nonumber\\
&&
+ C |D^k Rm|\cdot 
\Bigg[ 
\sum_{i=0}^k |D^i Rm|\cdot |D^{k-i}Rm| + \sum_{i=0}^k |D^i T|\cdot |D^{k+1 -i}Rm| 
\nonumber\\
&&
+ |D^{k+2}(T\ast T)| + \sum_{i=0}^{k-1} |D^i Rm|\cdot |D^{k-i}(T\ast T)|
\nonumber\\
&&
+ |D^{k+1}(T\ast Rm)|
+ |D^{k+1} (\bar T\ast Rm)|
 + \sum_{i=0}^{k-1} |D^i Rm|\cdot |D^{k-1-i} (\bar T\ast Rm)| 
\nonumber\\
&&
+|D^k (Rm\ast Rm)| + |D^{k+1} (T\ast \Psi)| + |D^{k}(\Psi\ast T\ast T) | + |D^k(\nabla\Psi\ast T)|\nonumber\\
&&
+ \sum_{i=1}^k |D^{k-i} H |\cdot |D^i \left({1\over 2\|\Omega\|_{\o}}\right)|
+ \sum_{i=1}^k |D^{k-i} Rm|\cdot |D^i (\p_t g)|
\Bigg]
\Bigg\}
\nonumber\\
&&
+ {C \over 2\|\Omega\|_{\o}} |D^k Rm|^2 \cdot |\Psi|
\eea

We estimate the terms on right hand side one by one. Recall that we have
\bea
|D^j Rm| &\leq & {C\, A\over t^{j/2}}, \ \ \ \  0\leq j \leq k-1\\
|D^{j+1} T| &\leq & {C\, A\over t^{j/2}}, \ \ \ \  0\leq j \leq k-1\\
|T|^2 &\leq& C\, A;
\eea
and the unknown terms are $|D^{k+1} Rm|, |D^k Rm|, |D^{k+2}T|$ and $|D^{k+1} T|$.

\

$\bullet$ Estimate for $| D^{k +1} Rm| \cdot \sum_{i=0}^{k-1} |D^i Rm|\cdot |D^{k -1-i} Rm|:$
\bea
| D^{k +1} Rm| \cdot \sum_{i=0}^{k-1} |D^i Rm|\cdot |D^{k -1-i} Rm| 
&\leq & | D^{k +1} Rm| \cdot \sum_{i=0}^{k-1}{CA\over t^{i/2}}\cdot {CA\over t^{(k-1-i)/2}}
\\\nonumber
&\leq &  | D^{k +1} Rm| \cdot CA^2 \, t^{-{k-1\over 2}}
\\\nonumber
&\leq & \theta  | D^{k +1} Rm|^2  + C(\theta) \, A^3\, t^{-{k}}
\eea
where $\theta$ is a small positive number such that $C_1\theta < {1\over 4}$ . To obtain the last inequality, we used Cauchy-Schwarz inequality and the fact that $A\,t<1$.

\medskip

$\bullet$ Estimate for $\sum_{i=0}^{k-1} |D^i Rm|^2\cdot |D^{k-1-i} Rm|^2:$
\bea
\sum_{i=0}^{k-1} |D^i Rm|^2\cdot |D^{k-1-i} Rm|^2 
&\leq &
\sum_{i=0}^{k-1} \left({CA\over t^{i/2}}\right)^2 \cdot \left( { CA\over t^{(k-1-i)/2}}\right)^2
\\\nonumber
&\leq & CA^4\, t^{-(k-1)} \leq CA^3 \, t^{-k}
\eea

\medskip

$\bullet$ Estimate for $\sum_{i=0}^k |D^i Rm|\cdot |D^{k-i}Rm| :$
\bea
\sum_{i=0}^k |D^i Rm|\cdot |D^{k-i}Rm| 
&=& 2 |D^k Rm|\cdot |Rm| + \sum_{i=1}^{k-1} |D^i Rm|\cdot |D^{k-i}Rm| 
\\\nonumber
&\leq& 
CA \, |D^k Rm| + CA^2 \, t^{-{k\over 2}}
\eea

\medskip

$\bullet$ Estimate for $\sum_{i=0}^k |D^i T|\cdot |D^{k+1 -i}Rm|:$
\bea
\sum_{i=0}^k |D^i T|\cdot |D^{k+1 -i}Rm|
&=& 
|T| \cdot |D^{k+1} Rm| + |DT| \cdot |D^k Rm| + \sum_{i=2}^k |D^i T|\cdot |D^{k+1 -i}Rm|
\nonumber\\
&\leq &
CA^{1\over 2}\, |D^{k+1} Rm| + CA \,  |D^k Rm|  + CA^2 \, t^{-{k\over 2}}
\eea

\medskip

$\bullet$ Estimate for $ |D^{k+2}(T\ast T)| :$
\bea
 |D^{k+2}(T\ast T)| 
 &\leq&
 \sum_{i=0}^{k+2} |D^i T| \cdot |D^{k+2-i} T|
 \\\nonumber
 &=&
 2|T| \cdot | D^{k+2}T| + 2 |DT|\cdot |D^{k+1}T| + \sum_{i=2}^{k} |D^i T|\cdot |D^{k+2-i} T|
  \\\nonumber
 &\leq & 
 CA^{1\over 2} \, | D^{k+2}T| + CA \, |D^{k+1} T| + CA^2 \, t^{-{k\over 2}}
\eea

\medskip

$\bullet$ Estimate for $\sum_{i=0}^{k-1}|D^i Rm|\cdot |D^{k-i}(T\ast T)|:$
\bea
&& \sum_{i=0}^{k-1} |D^i Rm|\cdot |D^{k-i}(T\ast T)|
\\\nonumber
&\leq &
2\sum_{i=0}^{k-1}  |D^i Rm|\cdot |T| \cdot |D^{k-i} T| 
+  \sum_{i=0}^{k-1} \sum_{j=1}^{k-i} |D^i Rm|\cdot |D^{j} T|\cdot |D^{k-i-j} T|
\\\nonumber
&\leq &
CA^2 \, t^{-{k\over 2}} 
\eea

\medskip

$\bullet$ Estimate for $|D^{k+1}(T\ast Rm)|:$
\bea\label{TRm}
|D^{k+1}(T\ast Rm)|
&\leq & |T| \cdot |D^{k+1} Rm| + |DT|\cdot |D^{k}Rm| + |D^{k+1} T|\cdot |Rm|
\\\nonumber
&& 
+\sum_{i=2}^{k} |D^{i} T|\cdot |D^{k+1 -i} Rm|
\\\nonumber
&\leq &
CA^{1\over 2}\, |D^{k+1} Rm| + CA\, |D^{k}Rm| + CA\,  |D^{k+1} T| + CA^2 \, t^{-{k\over 2}}
\eea

\medskip

$\bullet$ Estimate $ |D^{k+1} (\bar T\ast Rm)| :$
\bea
 |D^{k+1} (\bar T\ast Rm)| \leq CA^{1\over 2}\, |D^{k+1} Rm| + CA\, |D^{k}Rm| + CA\,  |D^{k+1} T| + CA^2 \, t^{-{k\over 2}}
\eea

\medskip

$\bullet$ Estimate for $ \sum_{i=0}^{k-1} |D^i Rm|\cdot |D^{k-1-i} (\bar T\ast Rm)|:$
\bea
&& \sum_{i=0}^{k-1} |D^i Rm|\cdot |D^{k-1-i} (\bar T\ast Rm)|\\\nonumber
 &\leq &
 \sum_{i=0}^{k-1}  |D^i Rm|\cdot |T| \cdot |D^{k-1-i} Rm| 
+  \sum_{i=0}^{k-1} \sum_{j=1}^{k-1-i} |D^i Rm|\cdot |D^{j} T|\cdot |D^{k-1-i-j} Rm|
 \\\nonumber
 &\leq &
 CA^2 \, t^{- {k\over 2}} 
\eea

\medskip

$\bullet$ Estimate for $|D^k (Rm\ast Rm)|:$
\bea
|D^k (Rm\ast Rm)| 
&\leq & 
2|Rm|  \cdot |D^k Rm| + \sum_{i=1}^{k-1} |D^i Rm|\cdot |D^{k-i} Rm|
\\\nonumber
&\leq & 
C A\, |D^k Rm| + CA^2 \, t^{- {k\over 2}} 
\eea

\medskip

$\bullet$ Estimate for $|D^{k+1} (T\ast \Psi)|:$

Recall that $\Psi_{\bar p q} = - \tilde R_{\bar p q} + g^{s\bar r} \, g^{m \bar n} \, T_{\bar n s q } \, \bar T_{m \bar r \bar p}$, we have
\bea
|D^{k+1} (\Psi\ast T) | 
&\leq& 
|D^{k+1} (Rm\ast T)| + |D^{k+1} ( T\ast T\ast T)|
\eea
The first term is the same as (\ref{TRm}). We only need to estimate the second term.
\bea
 |D^{k+1} ( T\ast T\ast T)| 
 &\leq &
 | D^{k+1} T |\cdot | T|_{\o}^2 + \sum_{p+q = k+1; p, q>0} |D^p T |\cdot |D^q T|\cdot |T|
 \\\nonumber
 &&
 + \sum_{p+q + r= k+1; p, q, r>0} | D^p T|\cdot |D^q T| \cdot |D^r T|
\\\nonumber
&\leq &
CA \, |D^{k+1} T| + CA^{5\over2} \, t^{- {(k-1)\over 2}} + C A^3 \, t^{- {(k-2)\over 2}}
\\\nonumber
&\leq &
CA \, |D^{k+1} T|  + CA^2\, t^{-{k\over 2}}
\eea
It follows that
\bea\label{PhiT1}
|D^{k+1} (\Psi\ast T) | 
&\leq& 
CA^{1\over 2}\, |D^{k+1}Rm| + CA\, (|D^k Rm| + |D^{k+1}T|) + C A^2 \, t^{-{k\over 2}}
\eea

\medskip

$\bullet$ Estimate for $|D^{k}(\Psi\ast T\ast T) |:$
\bea\label{PhiTT1}
|D^k(\Psi \ast T\ast T)|
&\leq &
|D^k ( Rm \ast T\ast T)| + |D^k ( T\ast T \ast T\ast T)|
\eea
We use the same trick as above to estimate these two terms. For the first term, we have
\bea
|D^k ( Rm \ast T\ast T)| 
&\leq &
|D^k Rm|\cdot |T |_{\o}^2 + \sum_{p+q=k; q>0} |D^p Rm| \cdot |D^q T| \cdot |T|
\\\nonumber
&&
+ \sum_{p+q+r=k; q, r>0} |D^p Rm| \cdot |D^q T |\cdot |D^r T|
\\\nonumber
&\leq &
CA\, |D^k Rm| + C A^{5\over 2} \, t^{-{k-1\over 2}} + CA^3 \, t^{-{k-2\over 2}}
\\\nonumber
&\leq &
CA\, |D^k Rm| + C A^{2} \, t^{-{k\over 2}}
\eea
For the second term, we have
\bea
 |D^k ( T\ast T \ast T\ast T)|
 &\leq &
 4|D^k T| \cdot |T|^3 + \sum_{p+q=k; \, p, q>0} |D^p T| \cdot |D^q T|\cdot | T|_{\o}^2
 \\\nonumber
 && 
 + \sum_{p+q + r=k; \, p, q, r>0} |D^p T|\cdot |D^q T|\cdot |D^r T|\cdot | T|_{\o}
  \\\nonumber
&&
   + \sum_{p+q + r+s=k;\, p, q, r, s>0} |D^p T|\cdot |D^q T|\cdot |D^r T|\cdot |D^s T| 
 \\\nonumber
 &\leq & 
 CA^{5\over 2} \, t^{- {k-1\over 2}} + CA^3 \, t^{-{k-2\over 2}} + CA^{7\over 2}\, t^{-{k-3\over 2}} + CA^4 \, t^{-{k-4\over 2}}
 \\\nonumber
 &\leq & CA^2 \, t^{-{k\over 2}}
\eea
Thus, we have
\bea\label{PhiTT1}
|D^k(\Psi \ast T\ast T)|
&\leq &
CA\, |D^k Rm| + C A^{2} \, t^{-{k\over 2}}.
\eea

\medskip

$\bullet$ Estimate for $|D^k(\nabla\Psi\ast T)|:$
\bea
|D^k (\nabla\Psi\ast T)| 
&\leq&
|D^k(\nabla Rm \ast T)| + |D^k (\nabla(T\ast T)\ast T)|
\\\nonumber
&\leq &
|D^{k+1}Rm|\cdot |T| + |D^k Rm|\cdot |DT| + \sum_{i=2}^k |D^{k+1-i} Rm|\cdot |D^i T|
\\\nonumber
&&
+  |D^{k+1} (T\ast T)| \cdot |T| + \sum_{i=1}^k |D^{k+1 - i } (T\ast T)| \cdot |D^i T|
\\\nonumber
&\leq & CA^{1\over 2}\, |D^{k+1}Rm| + CA\, |D^k Rm| + CA\, |D^{k+1} T|+ CA^2\, t^{-{k\over 2}} 
\eea

\medskip

$\bullet$ Estimate for $ \sum_{i=1}^k |D^{k-i} H |\cdot |D^i \left({1\over 2\|\Omega\|_{\o}}\right)|:$

Recall that 
\bea\label{Hformula}
H&=& {1\over 2}\Delta_{\R} Rm + \nabla\bar\nabla(T\ast \bar T) + \bar\nabla(T\ast Rm) + \nabla(\bar T \ast Rm)\\
  \nonumber
  &&
  + Rm\ast Rm + (\bar\nabla T - \bar T \ast T )\ast \Psi + \bar T \ast \nabla \Psi + T\ast \bar\nabla \Psi
\eea
and we also compute, for any $m$,
\bea\label{derivativeO1}
\nabla^m \left({1 \over 2\|\Omega\|_{\o}}\right) 
&=&
 \nabla^{m-1} \nabla \left({1 \over 2\|\Omega\|_{\o}}\right) = - \nabla^{m-1} \left({1 \over 2\|\Omega\|_{\o}} T\right)
 \\\nonumber
 &=& 
 - \nabla^{m-1} \left({1 \over 2\|\Omega\|_{\o}}\right) \ast T  -  {1 \over 2\|\Omega\|_{\o}} \nabla^{m -1} T
 \\\nonumber
 &=& {1 \over 2\|\Omega\|_{\o}} \sum_{j=1}^{m} \nabla^{m-j} T \ast T^{j-1}
\eea
where $T^{j-1} = T \ast T \ast \cdots \ast T$ with $(j-1)$ factors.
Again keep in mind that the unknown terms are $|D^{k+1}Rm|$, $|D^{k}Rm|$,$ |D^{k+2}T|$ and $|D^{k+1}T|$.
Notice that these terms only appear for $i=1, 2$ in the summation. 
\bea
\sum_{i=1}^k |D^{k-i} H |\cdot |D^i \left({1\over 2\|\Omega\|_{\o}}\right)|
&=& 
|D^{k-1} H |\cdot  |D \left({1\over 2\|\Omega\|_{\o}}\right)| + |D^{k-2} H |\cdot  |D^2 \left({1\over 2\|\Omega\|_{\o}}\right)| 
\nonumber\\
&&
+ \sum_{i=3}^k |D^{k-i} H |\cdot |D^i \left({1\over 2\|\Omega\|_{\o}}\right)|
\eea
Using (\ref{Hformula}) and (\ref{derivativeO1}), we can estimate the terms on the right hand side one by one and obtain
\bea
&&\sum_{i=1}^k |D^{k-i} H |\cdot |D^i \left({1\over 2\|\Omega\|_{\o}}\right)|
\\\nonumber
&\leq &
 CA^{1\over 2} |D^{k+1} Rm| + CA\, (|D^k Rm|+ |D^{k+1} T|) + CA^2 \, t^{-{k\over 2}}
\eea

\medskip

$\bullet$ Estimate for $\sum_{i=1}^k |D^{k-i} Rm|\cdot |D^i (\p_t g)|:$
\bea
|D^{i}(\p_t g) |
&=&
|D^{i} \left({1 \over 2\|\Omega\|_{\o}}\Psi \right)| = \sum_{j=0}^{i} |D^j  \left({1 \over 2\|\Omega\|_{\o}} \right) |\cdot |D^{i-j} \Psi|
\eea
By the definition of $\Psi$ and the computation (\ref{derivativeO1}), we know that the only unknown term appeared in the summation is when $j=i=k$. Thus, we arrive the following estimate
\bea
\sum_{i=1}^k |D^{k-i} Rm|\cdot |D^i (\p_t g)| 
&\leq & CA\, |D^k Rm| + CA^2 \, t^{-{k\over 2}}
\eea

\medskip

$\bullet$ Estimate for the last term $|D^k Rm|^2 \cdot |\Psi|:$
\bea
|D^k Rm|^2 \cdot |\Psi| \leq CA\, |D^k Rm|^2. 
\eea

Finally, putting all the above estimates together, we obtain the lemma. 
Q.E.D.

\

Following the same strategy, we can also prove the following lemma on estimates for the derivatives of the torsion.

\medskip

\begin{lemma}\label{Dk+1T}
Under the same assumption as in Lemma \ref{DkRm}, we have
\bea
\p_t |D^{k+1} T|^2 
&\leq &
{1 \over 2\|\Omega\|_{\o}}\, 
\Big\{
{1\over 2}\Delta_{\R} |D^{k+1} T|^2  - {3\over 4}\, |D^{k+2} T|^2 
\\\nonumber
&&
+CA^{1\over 2} \left( |D^{k+2} T| + |D^{k+1} Rm|\right) \cdot |D^{k+1} T|
\nonumber\\
&&
+ CA\left( |D^{k+1} T| + |D^k Rm|\right) \cdot |D^{k+1}T| \nonumber\\
&&+ CA^2 \, t^{-{k\over 2}}\, |\nabla^{k+1}T| + C A^3 t^{-k}
\Big\}.\nonumber
\eea
\end{lemma}

\

Now we return to the proof of Theorem \ref{Higherorderestimate}:

\medskip

We first prove the estimate (\ref{CTestimate}) for the case $k=1$. To obtain the desired estimate, we apply the maximum principle to the function
\bea
G_1(z, t) = t \left( |DRm|^2 + |D^2 T|^2\right) + \Lambda\left(|Rm|^2 + |DT|^2\right)
\eea

Using Lemma \ref{DkRm} and Lemma \ref{Dk+1T} with $k=1$, we have
\bea
&&\p_t \left( |D Rm|^2 + |D^{2} T|^2 \right) \\\nonumber
&\leq & 
{1 \over 2\|\Omega\|_{\o}}\, 
\Big\{
{1\over 2}\Delta_{\R} \left( |D Rm|^2 + |D^{2} T|^2\right)  -{3\over 4}\,  \left( |D^{2} Rm|^2 + |D^{3}T|^2\right)
\\\nonumber
&&
+ CA^{1\over 2}\, \left( |D^{2} Rm| + |D^{3} T|\right) \cdot \left( |D Rm| + |D^{2} T|\right)
\\\nonumber
&&
+ CA\, \left(|D Rm| + |D^{2}T|\right)^2 + CA^2 \, t^{-{1\over 2}}\,  \left(|D Rm| + |D^{2}T|\right) + C A^3 \, t^{-1}
\Big\}
\\\nonumber
&\leq &
{1 \over 2\|\Omega\|_{\o}}\, 
\Big\{
{1\over 2}\Delta_{\R} \left( |D Rm|^2 + |D^{2} T|^2\right)  - {1\over 2} \, \left( |D^{2} Rm|^2 + |D^{3}T|^2\right)
\\\nonumber
&&
+ CA\, \left(|D Rm|^2 + |D^{2}T|^2\right) + CA^3 \, t^{-1}
\Big\}
\eea
where we used the Cauchy-Schwarz inequality in the last inequality.

Recall the evolution equation
\bea
\p_t (|D T|^2 + |Rm|^2) 
&\leq & 
{1 \over 2\|\Omega\|_{\o}}\,  
\Big\{
{1\over 2}\Delta_{\R} ( |D T|^2 + |Rm|^2 ) - {1\over 2} ( | D^2 T|^2 + |D Rm|^2) 
\nonumber\\
&&
+CA^{3\over 2} ( |D Rm| + |D^2 T|) + CA^3
\Big\}.
\eea
It follows that
\bea
\p_t G_1
&\leq& 
{1 \over 4\|\Omega\|_{\o}}\,  
\Big\{
\Delta_{\R} G_1 
-t\left( |D^2 Rm|^2 + |D^3 T|^2 \right)
- \Lambda \left( |D^2 T|^2 + |D Rm|^2\right)
\nonumber\\
&&
+ CA\, t \, ( |D Rm|^2 + |D^2 T|^2 ) +CA^3
\\\nonumber
&&
+ CA^{3\over 2}\, \Lambda (|DRm| + |D^2 T|) + CA^3\, \Lambda
\Big\}
+ (|DRm|^2 + |D^2 T|^2)
\eea
Again, using Cauchy-Schwarz inequality,
\bea
C\, A^{3\over 2}\, \Lambda (|D Rm| + |D^2 T|)  
\leq  C A^3 \, \Lambda + \Lambda (|D Rm|^2 + |D^2 T|^2).
\eea
Putting these estimates together, we have
\bea
\p_t G 
&\leq & 
{1 \over 4\|\Omega\|_{\o}}\,  
\Big\{
\Delta_{\R} G - t\left( |D^2 Rm|^2 + |D^3 T|^2 \right)
\\\nonumber
&&
+(\|\Omega\|_{\o}- \Lambda  + C A t) \left( |D^2 T|^2 + |D Rm|^2\right)
+CA^3
\Big\}
\eea
By $A t \leq 1$ and choosing $\Lambda$ large enough, 
\bea
\p_t G 
&\leq & 
{1 \over 4\|\Omega\|_{\o}}\,  
\Big\{
2\Delta_{\R} G + CA^3 \Lambda
\Big\}
\eea
We note that the choice of constant $\Lambda$ depends on the upper bound of $\|\Omega\|_{\o}$. However, with the assumption (\ref{assumptiononcurvature}), we can get the uniform $C^0$ bound of the metric depending on the uniform lower bound of $\|\Omega\|_{\o}$. Consequently, we obtain the upper bound of $\|\Omega\|_{\o}$, which also depends on the uniform lower bound of $\|\Omega\|_{\o}$.

To finish the proof for $k=1$, observing that when $t=0$,
\bea
G(0) = {\Lambda\over 2} (|D T|^2 + |Rm|^2) \leq C \Lambda A^2.
\eea
Thus, applying the maximum principle to the above inequality implies that
\bea
G(t) \leq C\Lambda A^2 + CA^3 \Lambda \, t \leq C A^2
\eea
It follows
\bea
|D Rm| + |D^2 T| \leq {CA\over t^{1/2}}.
\eea
This establishes the estimate (\ref{CTestimate}) when $k=1$. Next, we use induction on $k$ to prove the higher order estimates.

\

Using Lemma \ref{DkRm} and Lemma \ref{Dk+1T} again, we have
\bea
&&\p_t \left( |D^k Rm|^2 + |D^{k+1} T|^2 \right) \\\nonumber
&\leq & 
{1 \over 2\|\Omega\|_{\o}}\, 
\Big\{
{1\over 2}\Delta_{\R} \left( |D^k Rm|^2 + |D^{k+1} T|^2\right)  -{3\over 4}\,  \left( |D^{k+1} Rm|^2 + |D^{k+2}T|^2\right)
\\\nonumber
&&
+ CA^{1\over 2}\, \left( |D^{k+1} Rm| + |D^{k+2} T|\right) \cdot \left( |D^k Rm| + |D^{k+1} T|\right)
\\\nonumber
&&
+ CA\, \left(|D^k Rm| + |D^{k+1}T|\right)^2 + CA^2 \, t^{-{k\over 2}}\,  \left(|D^k Rm| + |D^{k+1}T|\right) + C A^3 \, t^{-k}
\Big\}
\\\nonumber
&\leq &
{1 \over 2\|\Omega\|_{\o}}\, 
\Big\{
{1\over 2}\Delta_{\R}\left( |D^k Rm|^2 + |D^{k+1} T|^2\right)  - {1\over 2} \, \left( |D^{k+1} Rm|^2 + |D^{k+2}T|^2\right)
\\\nonumber
&&
+ CA\, \left(|D^k Rm|^2 + |D^{k+1}T|^2\right) + CA^3 \, t^{-k}
\Big\}.
\eea

Denote 
\bea\label{evolution f_k}
f_j(z, t) =  |D^j Rm|^2 + |D^{j+1} T|^2.
\eea
Then,
\bea
\p_t f_k \leq {1 \over 4\|\Omega\|_{\o}}\, \left( \Delta_{\R} f_k - f_{k+1} + CA\, f_k + CA^3 \, t^{-k}\right).
\eea

Next, we apply the maximum principle to the test function
\bea\label{defG_k}
G_k(z, t) = t^k f_k + \, \sum_{i=1}^k\Lambda_i \, B^k_i \, t^{k-i} \, f_{k-i}
\eea
where $\Lambda_i \,(1\leq i \leq k)$ are large numbers to be determined and $B_i^k = {(k-1)!\over (k-i)!}$.
We note that, for $1\leq i <k$, we still have an inequality similar to (\ref{evolution f_k}) for $f_{k-i}$.
\bea\label{evolution f_k-i}
\p_t f_{k-i} &\leq& 
{1 \over 4\|\Omega\|_{\o}}\, \left( 2\Delta_{\R} f_{k-i} - f_{k-i+1} + CA\, f_{k-i} + CA^3 \, t^{-(k-i)}\right)
\nonumber\\
&\leq &{1 \over 4\|\Omega\|_{\o}}\, \left( 2\Delta_{\R} f_{k-i} - f_{k-i+1}+ CA^3 \, t^{-(k-i)}\right)
\eea
where we used the induction condition (\ref{assumption1}) for the term $f_{k-i}$ when $1\leq i<k$.
From (\ref{evolution f_k}) and (\ref{evolution f_k-i}), we deduce
\bea
\p_t G_k 
&=& 
k t^{k-1} \, f_k + t^k \p_t f_k  
+ \sum_{i=1}^{k-1}\Lambda_i \, B^k_i\, (k-i) \, t^{k-i-1} \, f_{k-i} + \sum_{i=1}^k\Lambda_i \, B^k_i \, t^{k-i} \,\p_t f_{k-i}
\\\nonumber
&=&
k t^{k-1} \, f_k + {1 \over 4\|\Omega\|_{\o}}\, t^k \left( 2\Delta_{\R} f_k - f_{k+1} + CA\, f_k + CA^3 \, t^{-k}\right)
+\sum_{i=1}^{k-1}\Lambda_i \, B^k_i\, (k-i) \, t^{k-i-1} \, f_{k-i}  
\\\nonumber
&&
+  {1 \over 4\|\Omega\|_{\o}}\, \sum_{i=1}^k\Lambda_i \, B^k_i \, t^{k-i} \, \left( 2\Delta_{\R} f_{k-i} - f_{k-i+1}+ CA^3 \, t^{-(k-i)}\right)
\\\nonumber
&=& {1 \over 4\|\Omega\|_{\o}}\, 2\Delta_{\R} G_k -  {1 \over 4\|\Omega\|_{\o}}\, t^k \, f_{k+1} + t^{k-1}f_k \left( k + {CA\, t\over 4|\Omega|_{\o}}\right) + {1 \over 4\|\Omega\|_{\o}}\, CA^3\left(1 +  \sum_{i=1}^k\Lambda_i \, B^k_i  \right)
\\\nonumber
&&
+\sum_{i=1}^{k-1}\Lambda_i \, B^k_i\, (k-i) \, t^{k-i-1} \, f_{k-i}   -  {1 \over 4\|\Omega\|_{\o}}\, \sum_{i=1}^k\Lambda_i \, B^k_i \, t^{k-i}\, f_{k-i+1} 
\\\nonumber
&\leq & 
{1 \over 4\|\Omega\|_{\o}}\, 2\Delta_{\R} G_k + t^{k-1}f_k \left( k + {CA\, t\over 4\|\Omega\|_{\o}} -{\Lambda_1 \, B_1^k\over 4\|\Omega\|_{\o}} \right) + {1 \over 4\|\Omega\|_{\o}}\,CA^3 
\\\nonumber
&&
+\sum_{i=1}^{k-1}\Lambda_i \, B^k_i\, (k-i) \, t^{k-i-1} \, f_{k-i}   -  {1 \over 4\|\Omega\|_{\o}}\, \sum_{i=2}^k\Lambda_i \, B^k_i \, t^{k-i}\, f_{k-i+1} 
\eea
We note that the last two terms can be re-written as
\bea
&& \sum_{i=1}^{k-1}\Lambda_i \, B^k_i\, (k-i) \, t^{k-i-1} \, f_{k-i}   -  {1 \over 4\|\Omega\|_{\o}}\, \sum_{i=2}^k\Lambda_i \, B^k_i \, t^{k-i}\, f_{k-i+1} 
\\\nonumber
&=& 
\sum_{i=1}^{k-1} \left(\Lambda_i \, B^k_i\, (k-i) - {1 \over 4\|\Omega\|_{\o}}\,\Lambda_{i+1} \, B^k_{i+1} \right) t^{k-i-1} \, f_{k-i}
\\\nonumber
&=&
\sum_{i=1}^{k-1}  \left(\Lambda_i- {1 \over 4\|\Omega\|_{\o}}\,\Lambda_{i+1}\right)\, B^k_{i+1}\, t^{k-i-1} \, f_{k-i}
\eea
Thus, we obtain
\bea
\p_t G_k 
&\leq& 
{1 \over 4\|\Omega\|_{\o}}\, 2\Delta_{\R} G_k + t^{k-1}f_k \left( k + {CA\, t\over 4\|\Omega\|_{\o}} -{\Lambda_1 \, B_1^k\over 4\|\Omega\|_{\o}} \right) + {1 \over 4\|\Omega\|_{\o}}\,CA^3 
\\\nonumber
&&
+\sum_{i=1}^{k-1}  \left(\Lambda_i- {1 \over 4\|\Omega\|_{\o}}\,\Lambda_{i+1}\right)\, B^k_{i+1}\, t^{k-i-1} \, f_{k-i}
\eea
Choosing $\Lambda_1$ large enough and $\Lambda_i \leq {1 \over 4\|\Omega\|_{\o}}\,\Lambda_{i+1}$ for $1\leq i \leq k-1$, we have
\bea
\p_t G_k 
&\leq& 
{1 \over 4\|\Omega\|_{\o}}\, (2\Delta_{\R} G_k + CA^3)
\eea
Note that 
\bea
\max_{z\in M} G(z, 0) = \Lambda_k \, B_k^k f_0 = {(k-1)!\over 2} \, \Lambda_k \, (|Rm|^2 + |D T|^2) \leq C A^2
\eea
Applying the maximum principle to the inequality satisfied by $G_k$, we have
\bea
 \max_{z\in M} G(z, t) \leq CA^2 + C A^3 \, t \leq CA^2.
\eea
Finally, we get
\bea
|D^k Rm| + |D^{k+1} T| \leq CA\, t^{-{k\over 2}}.
\eea
The proof of Theorem \ref{Higherorderestimate} is complete. Q.E.D.

\

\subsection{Doubling estimates for the curvature and torsion}

Let 
\bea
f(z, t) = |D T|_{\o}^2 + |Rm|_{\o}^2 + |T|^4_{\o}
\eea
and denote $f(t) = \max_{z\in M} f(z, t)$. We can derive a doubling-time estimate for $f(t)$, which roughly says that $f(t)$ cannot blow up quickly. 

\begin{proposition}
There is a constant $C$ depending on a lower bound for $\|\Omega\|_\o$ such that 
\bea
 \max_{M} \left(|D T|^2 + |Rm|^2 + |T|^4\right)(t) \leq 4  \max_{M} \left(|D T|^2 + |Rm|^2 + |T|^4\right)(0)
\eea
\end{proposition}
for all $t \in [0, {1\over 4C f^{1\over 2}(0)}]$.

\medskip
\noindent
{\it Proof.}
The proof is standard and we apply the maximum principle to $f(z, t)$. 
Recall the evolution equations, by taking $k=0$ in (\ref{evolutionDkRm}),
\bea
\p_t |Rm|^2 
&\leq&
{1 \over 2\|\Omega\|_{\o}}\,  
\Big\{
{1\over 2}\Delta_{\R} |Rm|^2 -  |DRm|^2
+ C |D^2 T|\cdot |Rm|\cdot |T|
\\\nonumber
&&
+ C |D Rm| \cdot |Rm| \cdot |T| 
+ C|D T|^2 \cdot |Rm| + C |D T|\cdot |Rm|^2
\\\nonumber
&&
+ C|Rm|^3 + C |D T|\cdot |Rm|\cdot |T|^2 + C|Rm|^2 \cdot |T|^2 + C|Rm|\cdot |T|^4
\Big\}
\eea
We apply the Young's inequalities and get
\bea
\p_t|Rm|^2 
&\leq &
{1 \over 2\|\Omega\|_{\o}}\,  
\Big\{
{1\over 2}\Delta_{\R} |Rm|^2 -  {1 \over 2} |DRm|^2 + {1 \over 2} |D^2 T|^2 + C\left( |D T|^3 + |Rm|^3 + |T|^6 \right) 
\Big\}\nonumber\\
\eea
Similarly, considering the evolution equation for $|D T|^2$ and $|T|^2$, we can derive
\bea
\p_t |\nabla T|^2 
&\leq &
{1 \over 2\|\Omega\|_{\o}}\,   
\Big\{
{1\over 2}\Delta_{\R} |DT|^2 - {1 \over 2} |D^2 T|^2 + {1 \over 2} |D Rm|^2 + C\left( |D T|^3 + |Rm|^3 + |T|^6 \right) 
\Big\}\nonumber\\
\eea
and
\bea
\p_t |T|^4  
&\leq & {1 \over 2\|\Omega\|_{\o}}\, 
\Big\{
{1\over 2}\Delta_{\R} |T|^4 + C\left( |D T|^3 + |Rm|^3 + |T|^6 \right) 
\Big\}
\eea
Putting the above evolution equations together, we have
\bea
\p_t f(z, t)
&\leq &
{1 \over 2\|\Omega\|_{\o}}\, 
\Big\{
{1\over 2}\Delta_{\R} f +  C\left( |DT|^3 + |Rm|^3 + |T|^6 \right) 
\Big\}
\\\nonumber
&\leq & 
{1 \over 2\|\Omega\|_{\o}}\, \left( {1\over 2}\Delta_{\R} f + C f^{3\over 2}\right)
\eea
Finally, by the maximum principle, we have
\bea
 \p_t f(t) \leq {C \over 2\|\Omega\|_{\o}}\, f^{3\over 2}
\eea
which implies that
\bea
f(t) \leq {f(0) \over \left(1- 2C\, f^{1\over 2}(0) \, t\right)^2 }
\eea
Thus, as long as the flow exists and $t\leq {1- {1\over A} \over 2C f^{1\over 2}(0)}$, we have $f(t) \leq A^2 f(0)$. Q.E.D.

\

\subsection{A criterion for the long-time existence of the flow}

We can give now the proof of Theorem \ref{longtimeexistence}. We begin by observing that, under the given hypotheses, the metrics $\o(t)$ are uniformly equivalent for $t\in (T-\delta,T)$. Our goal is to show that the metrics are uniformly bounded in $C^\infty$ for some interval $t\in (T-\delta,T)$. This would imply the existence of the limit $\o(T)$ of a subsequence $\o(t_j)$ with $t_j\to T$. By the short-time existence theorem for the Anomaly flow proved in \cite{PPZ1}, it follows that the flow extends to $[0,T+\epsilon)$ for some $\epsilon>0$.

\subsubsection{$C^1$ bounds for the metric}

We need to establish the $C^\infty$ convergence of (subsequence of) the metrics $g_{\bar kj}(t)$ as $t\to T$. We have already noted the $C^0$ uniform boundedness of $g_{\bar kj}(t)$. In this section, we establish the $C^1$ bounds.
For this, we fix a reference metric $\hat g_{\bar kj}$ and introduce the relative endomorphism
\bea
h^j{}_m(t)=\hat g^{j\bar p}g_{\bar pm}(t).
\eea
The uniform $C^0$ bound of $g_{\bar kj}(t)$ is equivalent to the $C^0$ bound of $h(t)$. We need to estimate the derivatives of $h(t)$. For this, recall the curvature relation between two different metrics $g_{\bar kj}(t)$ and $\hat g_{\bar kj}$,
\bea
R_{\bar kj}{}^p{}_m=\hat R_{\bar kj}^p{}_m
-
\p_{\bar k}(h^p{}_q \hat \na_j h^p{}_m)
\eea
where $\hat \na$ denotes the covariant derivative with respect to $\hat g_{\bar kj}$. This relation can be viewed as a second order PDE in $h$, with bounded right hand sides because the curvature $R_{\bar kj}{}^p{}_m$ is assumed to be bounded, and which is uniformly elliptic because the metrics $g_{\bar kj}(t)$ are uniformly equivalent (and hence the relative endomorphisms $h(t)$ are uniformly bounded away from $0$ and $\infty$). It follows that
\bea
\|h\|_{C^{1,\alpha}}\leq C.
\eea

\subsubsection{$C^k$ bounds for the metric}
We will use the notation $G_k$ for the summation of norms squared of all combinations of $\hat\na^m \overline{\hat \na^{\ell}}$ acting on $g$ such that $m+\ell=k$. For example,
\be
G_2 = |\hat\na \hat\na g|^2 + |\hat\na \overline{\hat\na} g|^2+ |\overline{\hat\na} \overline{\hat\na} g|^2.
\ee
We introduce the tensor
\be
\Theta^k{}_{ij} = - g^{k \bar{\ell}} \hat\na_i g_{\bar{\ell} j},
\ee
which is the difference of the background connection and the evolving connection: $\Theta = \Gamma_0 -\Gamma$. We will use the notation $S_k$ for the summation of norms squared of all combinations of $\na^m \overline{\na^\ell}$ acting on $\Theta$ such that $m+\ell=k$. For example,
\be
S_2 = |\na \na \Theta|^2 + |\na \overline{\na} \Theta|^2+ |\overline{\na} \overline{\na} \Theta|^2.
\ee
Our evolution equation is
\be
\p_t g_{\bar{p} q} = {1\over 2 \| \Omega \|_\o}\, \Psi_{\bar p q} ,
\ee
where $\Psi_{\bar p q} = - \tilde R_{\bar p q} + g^{\alpha\bar\beta}g^{s\bar r}T_{\bar\beta sq}\bar T_{\alpha\bar r\bar p}$.

\begin{proposition}
Suppose all covariant derivatives of curvature and torsion of $g(t)$ with respect to the evolving connection $\nabla$ are bounded on $[0,T)$. Then all covariant derivatives of ${\Phi_{\bar p q} \over 2 \| \Omega \|_\o}$ with respect to the evolving connection $\nabla$ are bounded on $[0,T)$.
\end{proposition}
{\it Proof:} Compute
\be
\na^m \overline{\na}^\ell \bigg( {\Psi_{\bar p q} \over 2 \| \Omega \|_\o} \bigg) = {1 \over 2} \sum_{i \leq m}\sum_{ j \leq \ell} \na^i \overline{\na}^j \bigg({1 \over \| \Omega \|_\o} \bigg)\na^{m-i} \overline{\na}^{\ell-j} \Psi_{\bar p q}.
\ee
We have
\bea
\na^i \overline{\na}^j \bigg({1 \over \| \Omega \|_\o} \bigg) &=& - \na^i \overline{\na}^{j-1} \left({\overline{T} \over \| \Omega \|_\o}\right) \nonumber\\
&=&  {1 \over \| \Omega \|_\o} \sum \na^{i_1} \overline{\na}^{i_2} T^{i_3} \ast \na^{i_4} \overline{\na}^{i_5} \overline{T}^{i_6} \ast T^{i_7}\ast \overline{T}^{i_8} .
\eea
Since $\Psi$ is written in terms of curvature and torsion, and $\| \Omega \|_\o$ has a lower bound, the proposition follows. Q.E.D.
\begin{proposition}
Suppose all covariant derivatives of curvature and torsion of $g(t)$ with respect to the evolving connection $\nabla$ are bounded on $[0,T)$. If $G_i \leq C$ and $S_{i-1} \leq C$ for all non-negative integers $i \leq k$, then $G_{k+1} \leq C$ and $S_k \leq C$ on $[0,T)$.
\end{proposition}
{\it Proof:} By the previous proposition, all covariant derivatives of ${\Psi_{\bar p q} \over 2 \| \Omega \|_\o}$ with respect to the evolving connection $\nabla$ are bounded on $[0,T)$. Let $m+\ell =k+1$, and compute
\bea
\hat\na^m \overline{\hat\na^\ell} {\Psi_{\bar p q} \over 2 \| \Omega \|_\o} &=& (\na+\Theta)^m (\overline{\na} + \overline{\Theta})^\ell {\Psi_{\bar p q} \over 2 \| \Omega \|_\o}
= \na^m \overline{\na}^{\ell-1} \bigg(\overline{\Theta} \,{\Psi_{\bar p q} \over 2 \| \Omega \|_\o} \bigg) + O(1) \nonumber\\
&=& \na^m \overline{\na}^{\ell-1} \overline{\Theta}\cdot {\Psi_{\bar p q} \over 2 \| \Omega \|_\o}  + O(1),
\eea
where $O(1)$ represents terms which involve evolving covariant derivatives of ${\Psi_{\bar p q} \over 2 \| \Omega \|_\o}$ and up to $(k-1)$th order evolving covariant derivatives of $\Theta$, which are bounded by assumption. If $\ell=0$, the right-hand side is replaced by $\na^{m-1} \Theta\cdot {\Psi_{\bar p q} \over 2 \| \Omega \|_\o}$. Next, we compute
\bea \label{dPsi}
\na^m \overline{\na}^{\ell-1} \overline{\Theta}^{\bar{k}}_{\bar{i} \bar{j}} &=& - g^{\ell \bar{k}} \na^m \overline{\na}^{\ell-1} \hat\na_{\bar{i}} g_{\bar{j} \ell} \nonumber\\
&=&  - g^{\ell \bar{k}} (\hat\na - \Theta)^m (\overline{\hat\na} - \overline{\Theta})^{\ell-1} \hat\na_{\bar{i}} g_{\bar{j} \ell} \nonumber\\
&=&  - g^{\ell \bar{k}} \hat\na^m \overline{\hat\na^{\ell-1}} \, \hat\na_{\bar{i}} g_{\bar{j} \ell} + O(1).
\eea
It follows that
\be
\bigg| \hat\na^m \overline{\hat\na^\ell} \, {\Psi_{\bar p q} \over 2 \| \Omega \|_\o} \bigg| \leq C\left(1 + |\hat\na^m\overline{\hat\na^{\ell}}\,g|\right).
\ee
By differentiating the evolution equation and using the above estimate, we have
\be
\p_t |\hat\na^m\overline{\hat\na^{\ell}} g|_{\hat g}^2 \leq C\left(1 + |\hat\na^m \overline{\hat\na^{\ell}}\, g|_{\hat g}^2\right),
\ee
hence $|\hat\na^m \overline{\hat\na^{\ell}}\, g|$ has exponential growth. This proves $G_{k+1} \leq C$. Then $S_k \leq C$ now follows from (\ref{dPsi}), since $\na^m \overline{\na}^{\ell} \Theta= \overline{\overline{\na}^{m}\na^\ell  \overline{\Theta}}$
and we can exchange evolving covariant derivatives up to bounded terms. Q.E.D.

\bigskip

\par
By the $C^1$ bound on the metric, we have $G_1 \leq C$. We see that $S_0 = |\Theta| \leq C$ by definition of $\Theta$. Hence we can apply the previous proposition to deduce any estimate of the form
\be
|\hat\na^m \overline{\hat\na^\ell }\,g| \leq C.
\ee 
By differentiating the evolution equation with respect to time, we obtain
\be
\p_t^i \hat\na^m \overline{\hat\na^\ell} g = \hat\na^m \overline{\hat\na^\ell} \p_t^i \bigg( {\Psi_{\bar p q} \over 2 \| \Omega \|_\o} \bigg) .
\ee
Time derivatives of ${\Psi_{\bar p q} \over 2 \| \Omega \|_\o}$ can be expressed as time derivatives of connections, curvature and torsion, which in previous sections have been written as covariant derivatives of curvature and torsion. It follows that $\hat\na^m \overline{\hat\na^\ell} \p_t^i \bigg( {\Psi_{\bar p q} \over 2 \| \Omega \|_\o} \bigg)$ can be written in terms of evolving covariant derivatives of curvature and torsion, and hence is bounded. Therefore
\be
\left| \p_t^i \hat\na^m \overline{\hat\na^\ell}\, g \right| \leq C,
\ee
on $[0,T)$. Q.E.D.

\

\

\section{Appendix}
\begin{appendix}

\section{Conventions for differential forms}
\setcounter{equation}{0}

Let $\varphi$ be a $(p,q)$-form on the manifold $X$. We define its components
$\varphi_{\bar k_1\cdots\bar k_q j_1\cdots j_p}$ by
\bea
\varphi=
{1\over p!q!}
\sum \varphi_{\bar k_1\cdots\bar k_q j_1\cdots j_p}\,
dz^{j_p}\wedge\cdots\wedge dz^{j_1}\wedge
d\bar z^{k_q}\wedge\cdots\wedge d\bar z^{k_1}.
\eea
Although $\phi$ can be expressed in several ways under the above form, we reserve the notation $\varphi_{\bar k_1\cdots\bar k_q j_1\cdots j_p}$ for the uniquely defined coefficients $\varphi_{\bar k_1\cdots\bar k_q j_1\cdots j_p}$ which are anti-symmetric under permutation of any two of the barred indices or any two of the unbarred indices. 

To each Hermitian metric $g_{\bar kj}$ corresponds a Hermitian, positive $(1,1)$-form defined by
\bea
\o=ig_{\bar kj}\,dz^j\wedge d\bar z^k.
\eea
The Hermitian property $\overline{g_{\bar kj}}=g_{\bar jk}$ is then equivalent to the condition $\overline\o=\o$.

\section{Conventions for Chern unitary connections}
\setcounter{equation}{0}

Let $E\to X$ be a holomorphic vector bundle over a complex manifold $X$. Let $H_{\bar\alpha\beta}$ be a Hermitian metric on $E$. The Chern unitary connection is defined by
\bea
\na_{\bar k}V^\alpha=\p_{\bar k}V^\alpha, \qquad
\na_kV^\alpha=H^{\alpha\bar\gamma}\p_k(H_{\bar\gamma \beta}V^\beta)
\eea
for $V^\alpha$ any section of $E$.
Its curvature tensor is then defined by
\bea
[\na_j,\na_{\bar k}]V^\alpha=F_{\bar kj}{}^\alpha{}_\beta V^\beta.
\eea
Explicitly, we have
\bea
\na_kV^\alpha=\p_kV^\alpha+A_{k\beta}^\alpha V^\beta,
\qquad
A_{k\beta}^\alpha=H^{\alpha\bar\gamma}\p_kH_{\bar\gamma\beta}
\eea
and
\bea
F_{\bar kj}{}^\alpha{}_\beta=-\p_{\bar k}A_{j\beta}^\alpha
=
-\p_{\bar k}(H^{\alpha\bar\gamma}\p_jH_{\bar\gamma\beta}).
\eea
In particular, when $E=T^{1,0}(X)$, and $g_{\bar kj}$ is a Hermitian metric on $X$, we have the corresponding formulas
\bea
&&
\na_{\bar k}V^p=\p_{\bar k}V^p, \qquad
\na_kV^p=g^{p\bar m}\p_k(g_{\bar mq}V^q)
\nonumber\\
&&
[\na_j,\na_{\bar k}]V^p=R_{\bar kj}{}^p{}_q V^q
\nonumber\\
&&
R_{\bar kj}{}^p{}_q=-\p_{\bar k}A_{jq}^p
=
-\p_{\bar k}(g^{p\bar m}\p_jg_{\bar mq})
\eea
Our convention for the curvature form $Rm$ is 
\bea
Rm=R_{\bar kj}{}^p{}_q dz^j\wedge d\bar z^k.
\eea
It is the same as in \cite{FY1, FY2}, but it differs from that of \cite{S} by a factor of $i$.

\smallskip

When the metric on $X$ has torsion, the commutator identities $[\na_j,\na_k]$ for the Chern connections on any holomorphic vector bundle are given by
Hence for any tensor $A$, we have
\be
[\na_j,\na_k] A = T^\lambda{}_{jk} \na_\lambda A, \ \ [\na_{\bar{j}},\na_{\bar{k}}] A = \bar{T}^{\bar{\lambda}}{}_{\bar{j} \bar{k}} \na_{\bar{\lambda}} A.
\ee
Some useful examples are
\bea
\na_c \na_a \na_{\bar{b}} A_{\bar{i} j \bar{k} \ell} &=& \na_a \na_c \na_{\bar{b}} A_{\bar{i} j \bar{k} \ell} - T^\lambda{}_{ca} \na_\lambda \na_{\bar{b}} A_{\bar{i} j \bar{k} \ell} \nonumber\\
&=& \na_a \na_{\bar{b}} \na_c A_{\bar{i} j \bar{k} \ell} -T^\lambda{}_{ca} \na_\lambda \na_{\bar{b}} A_{\bar{i} j \bar{k} \ell} \nonumber\\
&&+ \na_a (R_{\bar{b} c  \bar{i}}{}^{\bar{\lambda}} A_{\bar{\lambda} j \bar{k} \ell}+ R_{\bar{b} c  \bar{k}}{}^{\bar{\lambda}} A_{\bar{i} j \bar{\lambda} \ell}- R_{\bar{b} c} {}^{\lambda}{}_j A_{\bar{i} \lambda \bar{k} \ell}- R_{\bar{b} c} {}^{\lambda}{}_\ell A_{\bar{i} j \bar{k} \lambda}  ).
\eea
and
\bea
\na_c \na_{\bar{d}} \na_a \na_{\bar{b}} A &=& \na_c \na_a \na_{\bar{d}} \na_{\bar{b}} A + \na (Rm \ast \bar{\na} A)\nonumber\\
&=& \na_a \na_c \na_{\bar{d}} \na_{\bar{b}} A + T \ast \na \bar{\na} \bar{\na} A + \na (Rm \ast \bar{\na} A) \nonumber\\
&=& \na_a \na_c \na_{\bar{b}} \na_{\bar{d}} A +  \na \na (\bar{T}\ast\bar{\na} A)+ T \ast \na \bar{\na} \bar{\na} A + \na (Rm \ast \bar{\na} A) \nonumber\\
&=& \na_a \na_{\bar{b}} \na_c \na_{\bar{d}} A +  \na \na (\bar{T}\ast\bar{\na} A)+ T \ast \na \bar{\na} \bar{\na} A + \na (Rm \ast \bar{\na} A). \nonumber\\
\eea
The general pattern is
\bea
\na^{(k)} \overline{\na}^{(\ell)} \na_a \na_{\bar{b}} A &=& \na_a \na_{\bar{b}} \na^{(k)} \overline{\na}^{(\ell)}  A + \sum_{\nu + \lambda= k} \sum_{\mu+\rho=\ell} \na^{(\nu)} \overline{\na}^{(\mu)} Rm \ast \na^{(\lambda)} \overline{\na}^{(\rho)} A \nonumber\\
&&+ \sum_{\nu + \lambda= k} \sum_{\mu+\rho=\ell+1} \na^{(\nu)} \overline{\na}^{(\mu)} T \ast \na^{(\lambda)} \overline{\na}^{(\rho)} A \nonumber\\
&&+ \sum_{\nu + \lambda= k+1} \sum_{\mu+\rho=\ell} \na^{(\nu)} \overline{\na}^{(\mu)} \overline{T} \ast \na^{(\lambda)} \overline{\na}^{(\rho)} A. 
\eea

\section{Identities for non-K\"ahler metrics}
\setcounter{equation}{0}

When the metric is not K\"ahler, the integration by parts formula becomes
\bea
\int_X \na_j V^j \,\o^n
=
\int_X (A_{pj}^p-A_{jp}^p)V^j\,\o^n=
\int_X g^{p\bar q}T_{\bar qpj}\,V^j\,\o^n.
\eea
It is convenient to introduce
\bea
T_j=g^{p\bar q}T_{\bar qpj}
\eea
so that the above equation becomes
\bea
\int_X \na_jV^j\,\o^n=\int_X  T_j V^j\,\o^n.
\eea

\subsection{The adjoints $\bar\p^\dagger$ and $\p^\dagger$ with torsion}

Since the signs are crucial, we work out in detail the operators $\bar\p^\dagger$ and $\p^\dagger$ on the space $\Lambda^{1,1}$ of $(1,1)$-forms.

\smallskip
Consider first the operator $\bar\p:\Lambda^{1,0}\to\Lambda^{1,1}$. Explicitly,
\bea
\bar\p(f_jdz^j)
=
\p_{\bar k}f_j d\bar z^k\wedge dz^j
=
-\p_{\bar k}f_j dz^j\wedge d\bar z^k
\eea
which means that
\bea
(\bar\p f)_{\bar kj}=-\p_{\bar k}f_j.
\eea
Let $\Phi=\Phi_{\bar pq}dz^q\wedge d\bar z^p$ be a $(1,1)$-form. The adjoint $\bar\p^\dagger$ is characterized by the equation
\bea
\<\bar\p f,\Phi\>=\<f,\bar\p^\dagger\Phi\>
\eea
which is equivalent to
\bea
\int_X(-\p_{\bar k}f_j)\overline{\Phi_{\bar pq}}g^{p\bar k}g^{j\bar q}{\o^n\over n!}
=
\int_Xf_j\overline{(\bar\p^\dagger \Phi)_q}g^{j\bar q}{\o^n\over n!}.
\eea
Integrating by parts, we find
\bea
(\bar\p^\dagger \Phi)_q
=
g^{k\bar p}(\na_k\Phi_{\bar pq}
-
T^j{}_{kj}\Phi_{\bar pq})
=
g^{k\bar p}(\na_k\Phi_{\bar pq}-T_k\Phi_{\bar pq}).
\eea
Similarly, we work out $\p^\dagger$. For $f=f_{\bar k}d\bar z^k$, we have
$\p f=\p_jf_{\bar k} dz^j\wedge d\bar z^k$, so that $(\p f)_{\bar kj}=\p_jf_{\bar k}$. Thus, the equation $\<\p f,\Phi\>=\<f,\p^\dagger\Phi\>$ becomes
\bea
\int_X \p_jf_{\bar k}
\overline{\Phi_{\bar pq}}g^{p\bar k}g^{j\bar q}{\o^n\over n!}
=
\int_Xf_{\bar k}
\overline{(\p^\dagger \Phi)_{\bar p}}g^{p\bar k}{\o^n\over n!}.
\eea
This results now into
\bea
(\p^\dagger\Phi)_{\bar q}
=
-g^{p\bar j}(\na_{\bar j}\Phi_{\bar qp}-\bar T_{\bar j}\Phi_{\bar qp}).
\eea

\subsection{Bianchi identities for non-K\"ahler metrics}

It is well-known that the Riemann curvature tensor of K\"ahler metrics satisfies the following important identities
\bea
&&
R_{\bar \ell m\bar kj}
=
R_{\bar k m\bar\ell j}=R_{\bar kj\bar\ell m}
\nonumber\\
&&
\na_qR_{\bar\ell m}{}^k{}_j=\na_m R_{\bar\ell q}{}^k{}_j,
\qquad
\na_{\bar p}R_{\bar\ell m}{}^k{}_j=\na_{\bar\ell}R_{\bar pm}{}^k{}_j.
\eea
For general Hermitian metrics, these identities become
\bea
\label{bianchi1}
R_{\bar\ell m\bar kj}&=&R_{\bar\ell j\bar km}+\na_{\bar\ell}T_{\bar k jm}
\nonumber\\
R_{\bar\ell m\bar kj}&=&R_{\bar km\bar\ell j}+\na_m\bar T_{j\bar k\bar\ell}
\eea
and
\bea
\label{bianchi2}
&&
\na_m R_{\bar kj}{}^p{}_q
=
\na_jR_{\bar km}{}^p{}_q+T^r{}_{jm}R_{\bar kr}{}^p{}_q,
\qquad
\na_m R_{\bar kj}{}_{\bar pq}
=
\na_jR_{\bar km}{}_{\bar pq}+T^r{}_{jm}R_{\bar kr\bar pq}
\nonumber\\
&&
\na_{\bar m}R_{\bar kj}{}^p{}_q
=
\na_{\bar k}R_{\bar mj}{}^p{}_q
+
\bar T^{\bar r}{}_{\bar k\bar m}R_{\bar rj}{}^p{}_q,
\qquad
\na_{\bar m}R_{\bar kj\bar pq}
=
\na_{\bar k}R_{\bar mj\bar pq}
+
\bar T^{\bar r}{}_{\bar k\bar m}R_{\bar rj\bar pq}
\nonumber\\
\eea
Observe that to interchange, say $m$ and $q$ in the second Bianchi identity for non-K\"ahler metrics, we have to use the first Bianchi identity and differentiate, resulting into
\bea
\na_m R_{\bar kj}{}_{\bar pq}
-
\na_qR_{\bar kj}{}_{\bar pm}
=
\na_q\na_{\bar k}T_{\bar pmj}
+
\na_m\na_{\bar k}T_{\bar pqj}
+
T^r{}_{qm}R_{\bar kr\bar pj}.
\eea
The occurrence of $D^2T$ on the right hand side is a source of potential difficulties, so it is desirable not to exchange this type of pairs of indices.

\end{appendix}

\bigskip
\noindent
{\bf Acknowledgements} The authors would like to thank the referees for a particularly careful reading of their paper, and in particular for pointing out several typos that could have been quite confusing to the reader.

\bigskip
Department of Mathematics, Columbia University, New York, NY 10027, USA

\smallskip

phong@math.columbia.edu

\bigskip
Department of Mathematics, Columbia University, New York, NY 10027, USA

\smallskip
picard@math.columbia.edu

\bigskip
Department of Mathematics, University of California, Irvine, CA 92697, USA

\smallskip
xiangwen@math.uci.edu

\end{document}